\documentclass[10pt]{amsart}
\usepackage{epsfig,url}
\usepackage{mathdots}
\usepackage{xcolor}
\usepackage{marginnote}
\usepackage{hyperref}
\usepackage{tikz}
\usepackage{graphicx}
\usepackage{float}
\usepackage[a4paper, margin=1.4in]{geometry}
\usetikzlibrary{decorations.markings, matrix}

\newtheorem{theorem}{Theorem}[section]

\newtheorem{conjecture}[theorem]{Conjecture}

\newtheorem{corollary}[theorem]{Corollary}
\newtheorem{proposition}[theorem]{Proposition}
\newtheorem{RHP problem}[theorem]{Riemann Hilbert Problem}

\theoremstyle{definition}

\theoremstyle{remark}
\newtheorem{remark}[theorem]{Remark}

\newcommand{\abs}[1]{\lvert#1\rvert}

\begin{document}

\title{Asymptotics and zeros of a special family of Jacobi polynomials}


\author[]{John Lopez Santander}
\address{Department of Mathematics,
Tulane University, New Orleans, LA 70118}
\email{jlopez12@tulane.edu}

\author[]{Kenneth D. T-R McLaughlin}
\address{Department of Mathematics,
Tulane University, New Orleans, LA 70118}
\email{kmclaughlin@tulane.edu}

\author[]{Victor H. Moll}
\address{Department of Mathematics,
Tulane University, New Orleans, LA 70118}
\email{vhm@tulane.edu}


\date{\today}


\begin{abstract}
 In this paper we study a family of non-classical Jacobi polynomials with varying parameters of the form $\alpha_n=n+1/2$ and $\beta_n=-n-1/2$. We obtain global asymptotics for these polynomials, and use this to establish results on the location their zeros. The analysis is based on the Riemann Hilbert formulation of Jacobi polynomials derived from the non-hermitian orthogonality  introduced by Kuijlaars, et al. This family of polynomials arise in the symbolic evaluation integrals in the work of Boros and Moll and corresponds to a limitting case, which is not considered in the works of Kuijlaars, et al. A remarkable feature in the analyisis is encountered when performing the local analysis of the RHP near the origin, where the local parametrix introduces a pole.
\end{abstract}

\maketitle

\newcommand\norm[1]{\left\lVert#1\right\rVert}
\newcommand{\ba}{\begin{eqnarray}}
\newcommand{\ea}{\end{eqnarray}}
\newcommand{\ift}{\int_{0}^{\infty}}
\newcommand{\nn}{\nonumber}
\newcommand{\no}{\noindent}
\newcommand{\lf}{\left\lfloor}
\newcommand{\rf}{\right\rfloor}
\newcommand{\realpart}{\mathop{\rm Re}\nolimits}
\newcommand{\imagpart}{\mathop{\rm Im}\nolimits}
\newcommand{\K}{\mathbf{K}}
\newcommand{\J}{\mathbf{J}}
\newcommand{\A}{\mathbf{A}}

\newcommand{\op}[1]{\ensuremath{\operatorname{#1}}}
\newcommand{\pFq}[5]{\ensuremath{{}_{#1}F_{#2} \left( \genfrac{}{}{0pt}{}{#3}
{#4} \bigg| {#5} \right)}}

\newcommand{\pFqcomma}{\mskip\pFqmuskip}

\newtheorem{Definition}{\bf Definition}[section]
\newtheorem{Thm}[Definition]{\bf Theorem}
\newtheorem{Example}[Definition]{\bf Example}
\newtheorem{Lem}[Definition]{\bf Lemma}
\newtheorem{Cor}[Definition]{\bf Corollary}
\newtheorem{Prop}[Definition]{\bf Proposition}
\numberwithin{equation}{section}

\section{Introduction}
\label{sec-introduction}

The classical Jacobi polynomials $P_{n}^{(\alpha, \beta)}(x)$, defined for $\alpha, \, \beta > -1$, can be constructed by applying 
the Gram-Schmidt orthogonalization procedure to the standard basis $\{ 1, \, x, \, x^{2}, \ldots, x^{n}, \ldots \}$ with inner product 
\begin{equation}
\label{innerprod-1}
\langle p(x), q(x) \rangle = \int_{-1}^{1} p(x) q(x) w(x;\alpha, \beta) \, dx 
\end{equation}
and weight function 
\begin{equation}
\label{weight-1}
w(x;\alpha, \beta) = (1-x)^{\alpha} (1+x)^{\beta}.
\end{equation}
\noindent
The restriction on the parameters $\alpha, \, \beta$ imposed above guarantee the convergence of the integral in \eqref{innerprod-1} for 
arbitrary polynomials $p, \, q$.

Expressions for the Jacobi polynomials include the explicit formula 
\begin{equation}
\label{jaco-explicit1}
P_{n}^{(\alpha,\beta)}(x) = 2^{-n} \sum_{k=0}^{n} \binom{n+\alpha}{n-k} \binom{n+\beta}{k} (x-1)^{k} (x+1)^{n-k},
\end{equation}
\noindent
 the Rodriguez formula 
\begin{equation}
P_{n}^{(\alpha,\beta)}(x) = \frac{1}{2^{n} n!} (x-1)^{-\alpha} (x+1)^{-\beta} \left( \frac{d}{dx} \right)^{n} 
\left[ (x-1)^{n+ \alpha} (x+1)^{n+\beta} \right],
\end{equation}
\noindent
as well as the hypergeometric representation 
\begin{equation}
P_{n}^{(\alpha,\beta)}(x) = \frac{(\alpha+1)_{n}}{n!} 
\pFq21{-n \quad \,\,\, n+ \alpha + \beta +1}{\alpha + 1}{ \frac{1-x}{2}},
\end{equation}
\noindent 
with the hypergeometric function $_{2}F_{1}$ given by 
\begin{equation}
\pFq21{a \quad b }{c}{\, x} = \sum_{k=0}^{\infty} \frac{(a)_{k} \, (b)_{k}}{(c)_{k}} \frac{x^{k}}{k!}
\end{equation}
\noindent
and the Pochhammer symbol by 
\begin{equation}
(u)_{k} = \begin{cases} u(u+1) \cdots (u+k-1) & \quad \textnormal{for} \,\, k =1, \, 2, \, \ldots \\
1 & \quad \textnormal{for} \,\, k = 0. 
\end{cases}
\end{equation}

The Jacobi polynomials are unique up to a scaling constant. The expressions given above satisfy 
\begin{equation}
P_{n}^{(\alpha, \beta)}(1) = \frac{(\alpha+1)_{n}}{n!}.
\end{equation}
\noindent
More information about them appears in \cite{andrews1999special}. 

\begin{figure}[htbp]
\centering
    \includegraphics[width=0.3\textwidth]{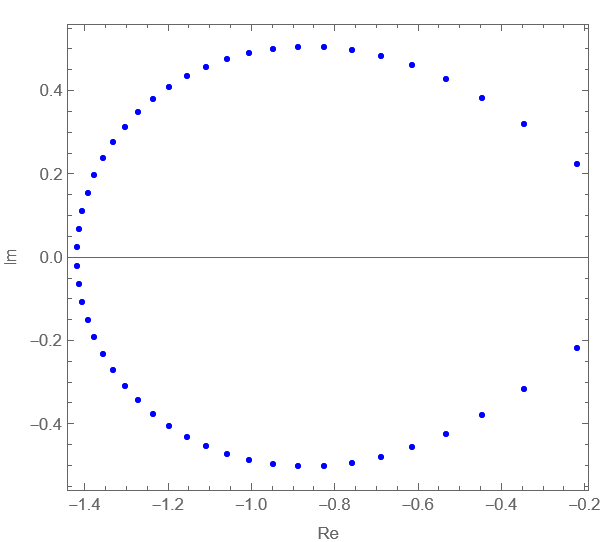}
    \caption{Zeros of $P_{50}(z)$.}
    \label{fig_zeros_p50}
\end{figure}

In the present work we consider the Jacobi polynomials 
with parameters $\alpha, \, \beta$ outside the classical range $\alpha, \, \beta > -1$. Specifically we define $P_m(z)$ by \eqref{jaco-explicit1}, with (varying) parameters $\alpha_m= m+ \tfrac{1}{2}$ and $\beta_m = -  m - \tfrac{1}{2}$: 
\begin{equation}
    \label{jaco-explicit2}
    P_m(z)=P_m^{\left(m+\frac{1}{2},-m-\frac{1}{2}\right)}(z).
\end{equation}
These polynomials appeared in the evaluation of the integral 
\begin{equation}
N_{0,4}(a;m) = \int_{0}^{\infty} \frac{dx}{(x^4+2ax^{2} + 1)^{m+1}},
\end{equation}
\noindent
given in \cite{boros1999integral, boros1999sequence} in the form 
\begin{equation}
N_{0,4}(a;m) = \frac{\pi}{2^{m+3/2} (a+1)^{m+\frac{1}{2}}}P_{m}(a).
\end{equation}
The following result regarding the zeros of these polynomials is proven in \cite{boros1999sequence}:
\begin{proposition}
\begin{enumerate}
    \item The polynomials $P_m(z)$ have no real zeros if $m$ even and a single real zero $x^*$, satisfying $x^*<-1$, if $n$ is odd.
    \item Let $\{z_k\}$, $1 \le k\le m$ be the sequence of zeros of $P_m(z)$. Then $|z_k + 1| < 1$.
\end{enumerate}    
\end{proposition}
The following conjecture was proposed in \cite{boros1999sequence}.
\begin{conjecture}\label{cojecture distance to -1}
The zeros $z_k$ of the polynomials $P_m(z)$ satisfy
\begin{equation}
   \lim_{m\rightarrow \infty} \max\{|z_k+1|: 1\le k\le m \}= 1.
\end{equation}
\end{conjecture}  
Moreover, as the degree of the polynomials grows to infinity, the zeros of the polynomials $P_m(z)$ seem to concentrate on a lemniscate. This phenomena has been studied for the partial sums of the exponential function in \cite{kriecherbauer2007locating}. Figure \ref{fig_zeros_p50} shows the zeros of $P_{50}(z)$.

The family $\{P_m(z)\}$ is not orthogonal with respect to a fixed weight $w(x)$ on the real line. This follows, for instance, from the fact that the zeros are not real. However, each polynomial $P_{m}(z)$ belongs to a different family $\{p_{m,n}\}_{n=0}^{\infty}$, orthogonal with respect to the weight function $w_{m}(z)=(1-z)^{m+\frac{1}{2}}(1+z)^{-m-\frac{1}{2}}$, where orthogonality is defined by integration on a contour in the complex plane, see Section \ref{sec-orthogonality} for more details.

The family of polynomials $\{P_m(z)\}$ is normalized by
\begin{equation}\label{Eq def pi_m}
\pi_m(z)=\frac{1}{\kappa_{m}} P_m^{\left(m+\frac{1}{2},-m-\frac{1}{2}\right)}(z),
\end{equation}
where the scaling constant $$\kappa_{m}=\frac{1}{2^m}\binom{2m}{m}$$ 
is chosen so that $\pi_m$ is monic.

The orthogonality relation in \eqref{Eq. orthogonality} and the associated Riemann Hilbert problem (see RHP \ref{RHP Y}), are used to obtain the global asymptotic behavior of $\{\pi_m(z)\}$ as $m\rightarrow\infty$. This is then used to derive the asymptotic behavior of the zeros of $\pi_m$.
Conjecture \ref{cojecture distance to -1} is established as the lemniscate behavior of the zeros. 

\begin{figure}[H]
    \centering
    \includegraphics[width=0.7\textwidth]{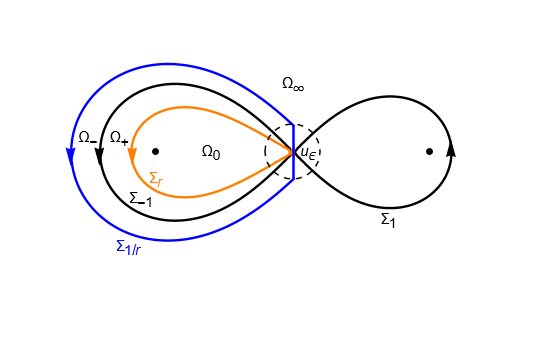} 
    \caption{Contour $\Sigma_U$.}
    \label{fig: Contour U}
\end{figure}

This analysis uses a decomposition of $\mathbb{C}$ into a number of regions given next, see also Figure \ref{fig: Contour U} for an illustration.

Consider the contour
\begin{equation}\label{Eq Contour Sigma}
   \Sigma=\Sigma_{-1}\cup\Sigma_1,
\end{equation}
where 
\begin{equation*}
   \Sigma_{-1}=\left\{|1 - z^2| = 1, \,\, \Re(z) \leq 0 \right\} \hspace{0,3cm}\text{and}\hspace{0.3cm} \Sigma_{-1}=\left\{  \left| 1 - \frac{3}{2} z^2 \right| = 1, \,\, \Re(z) \geq 0 \right\}.
\end{equation*}
Both $\Sigma_{-1}$ and $\Sigma_1$ are closed curves,  corresponding to the left and right halves of two distinct lemniscates.

Then, $\Sigma$ is extended to the contour $\Sigma_U$ shown in Figure \ref{fig: Contour U}. This is done by adding a closed, simple curve $\Sigma_r$ in the interior of $\Sigma_{-1}$. This curve emanates from the origin at a different angle than $\Sigma_{-1}$, stays at a finite distance from $\Sigma_{-1}$, and encircles $-1$. Additionally, we add a second curve, $\Sigma_{1/r}$, in the exterior of $\Sigma_{-1}$, consisting of a small segment of the imaginary axis, centered at the origin, and a simple curve joining the endpoints of the segment. This second curve remains at a finite distance from $\Sigma_{-1}$. All the contours are differentiable and oriented in the positive (counterclockwise) direction.

The contour $\Sigma_U$ creates the following regions in the complex plane: 
\begin{itemize}
    \item the region $\Omega_0$, inside $\Sigma_{r}$,
    \item the region $\Omega_+$, between $\Sigma_{r}$ and $\Sigma_{-1}$,
    \item  the region $\Omega_-$, between $\Sigma_{-1}$ and $\Sigma_{1/r}$, and,
    \item the region $\Omega_{\infty}$, outside $\Sigma_{1/r}$ except for the curve $\Sigma_{1}$.
\end{itemize}
Finally, we consider a small disk $u_{\varepsilon}$ centered at the origin and with the diameter being the segment of the imaginary axis in the curve $\Sigma_{1/r}$.

Let $D_{-\frac{1}{2}}(\,\cdot\,)$ be the parabolic cylinder function, defined in \cite{Bateman1953}[pag. 116], and let $\xi=\xi(z)$ the local map defined in \eqref{Eq def xi(z)}. The main results of this paper are given next:

\begin{theorem}\label{Thm Asymptotics}
     Let $\pi_m(z)$ be the monic Jacobi polynomial defined in \eqref{Eq def pi_m}. There exists a small neighbourhood $u_{\varepsilon}$ of the origin such that for $z\in \Omega_-\setminus u_{\varepsilon}$:
\begin{equation*}
\scalebox{1.1}{$
\pi_m(z)=\frac{(1+z)^{\frac{1}{2}}}{z^{\frac{1}{2}}} (1+z)^m\left(1-\frac{\sqrt{2}i}{(1-z^2)^m(1-z)^{\frac{1}{2}}(1+z)^{\frac{1}{2}}}+\mathcal{O}\left(\frac{1}{m}\right)\right)
$}.
\end{equation*}


\end{theorem}

For $z$ near the origin, let $ \xi=\lambda(z)$ be difined by the map in \eqref{Eq def xi(z)}, which has the local behavior $\xi\approx \sqrt{2m}z$. The asymptotics of the polynomials are described for a disk of fixed radius $R$ in the $\xi$ variable. 
\begin{theorem}\label{Thm Asymptotics local} Let $R>0$ be arbitrary. For  $|\xi|< R$ and $z=z(\xi)\in \Omega_-\cap u_{\varepsilon}\cap \mathbb{C}^+$,
\begin{equation*}
\begin{split}
\pi_m(z)& =\hspace{-0.05cm} \frac{(1+z)^m(1+z)^{\frac{1}{2}}\left(\xi\,e^{-\frac{5\pi i}{4}}\right)^{\frac{1}{2}} e^{\frac{i\xi^2}{4}}}{2\,z^{\frac{1}{2}}}\Bigg\{\left[1-\frac{1}{(1-z)^{\frac{1}{2}}(1+z)^{\frac{1}{2}}}\right]\\
&\hspace{2cm}\cdot \left[e^{\frac{\pi i}{4}}D_{-\frac{1}{2}}\left(e^{-\frac{3 \pi i}{4}}\xi\right)-i\sqrt{2}D_{-\frac{1}{2}}\left(e^{-\frac{5 \pi i}{4}}\xi\right)\right]+\mathcal{O}\left(\frac{1}{\sqrt{m}}\right)\Bigg\}.
\end{split} 
\end{equation*}
\end{theorem}

\begin{remark}
Throughout the paper, the standard notation $f_m(z) = \mathcal{O}\left(\frac{1}{m}\right)$ is used; there is a constant $C > 0$ (independent of $m$), and an integer $M > 0$ such that  
$$
\left| f_m(z) \right| \leq \frac{C}{m},
$$  
for all $m \geq M$, uniformly in $z$. 
\end{remark}

\begin{remark}
Similar expressions can be obtained for the other regions, $\Omega_0$, $\Omega_+$ and $\Omega_{\infty}$, both inside and outside $u_{\varepsilon}$. Theorem \ref{Thm Asymptotics} describes the region where the zeros of the polynomials $\pi_m$ are located. This is our primary application of Theorem \ref{Thm Asymptotics}. The interested reader will find in Section \ref{sec-RH Analysis} more information on the remaining regions.


\end{remark}

\begin{theorem}\label{Thm Zeros out}
Each zero of $\pi_m(z)$ in the region $\Omega_-\setminus u_{\varepsilon}$  approaches a zero of
$$q(m,z)=1-\frac{i\sqrt{2}}{(1-z^2)^m(1-z)^{\frac{1}{2}}(1+z)^{\frac{1}{2}}},\hspace{0,3cm} \text{as $m\rightarrow\infty$}.$$
More specifically, introduce the parametrization $1 - z^2 = \rho e^{i\theta}, \hspace{0.2cm} 0 \le \theta < 2\pi$. The zeros $z_k^*$ of $q(m,z)$ have $\rho =2^{\frac{1}{2m+1}}$ and can be enumerated according to \linebreak $\theta_k=\textnormal{Arg}(1-(z_k^*)^2)$, for $k=1, 2,\ldots,m$, where
    $$\left(m+\frac{1}{2}\right)\theta_k=\frac{\pi}{2}+2\pi k.$$
There exist an integer $M>0$ and a constant $C > 0$, such that if $m\ge M$, then for each zero $z_k^*$ of $q(m,z)$ outside $u_{\varepsilon}$, there is exactly one zero $z_k$ of $\pi_m(z)$, satisfying
$$|z_k^*-z_k|\le \frac{C}{m}.$$
\end{theorem}

This yields the following asymptotic result on the zeros of the polynomials $P_m(z)$.

\begin{corollary}\label{Corollary Asymp Zeros}
The zeros of the polynomials $P_m(z)$, outside the neighbourhood $u_{\varepsilon}$, approach the lemniscate $|1-z^2|=1$. More precisely, let $Z^{(m)} = \{z_k^{(m)} \in \mathbb{C} \setminus u_{\varepsilon} : \, z_k^{(m)} \text{ is a zero of } P_m(z)\}$, then
\begin{equation*}
  \max_{z_k^{(m)}\in Z^{(m)}}  \,\min\{|z_k^{(m)}-z|:z\in \mathbb{C},\,\, |1-z^2|=1\}\rightarrow 0,\hspace{1cm}\text{as $m\rightarrow\infty$}.
\end{equation*}  
\end{corollary}

This solves Conjecture \ref{cojecture distance to -1}.

\begin{remark}
This is consistent with the results of Driver and Möller \cite{driver2001zeros}, who used the hypergeometric representation of Jacobi polynomials to established results on the location of the zeros of non-classical 
 Jacobi polynomials in the same parameter regime. In their work, they applied connection formulae to transform the equation $P_m(z)= 0$, proving that the limit curve of the reciprocal of the zeros is a Cassini curve, see \cite[Thm. 4.1]{driver2001zeros}. However, their methods do not provide information about the argument of the roots, see \cite[page 86]{driver2001zeros}. A key advantage of our approach is that it describes how the zeros of $P_m(z)$ and $q(m, z)$ are asymptotically close. Since the zeros of the function $q(m, z)$ can be determined explicitly, our method yields information about both the modulus and the argument of the roots of $P_m(z)$ .
\end{remark}

In the next statement, the zeros of $\pi_m$ near to the origin are related to the zeros of a function involving the parabolic cylinder function. This function is expressed in terms of the variable $\xi=\xi(z)$ defined by the map $\lambda$ in \eqref{Eq def xi(z)}.

\begin{theorem}\label{Thm Zeros in}
Let $R>0$ be arbitrary. For each zero $\xi_k^*$ of 
$$q_{\ell}(\xi)=e^{\frac{\pi i}{4}}D_{-\frac{1}{2}}\left(e^{-\frac{3 \pi i}{4}}\xi\right)-i\sqrt{2}D_{-\frac{1}{2}}\left(e^{-\frac{5 \pi i}{4}}\xi\right),$$
with $|\xi_k^*|< R$ and $z(\xi_k^*)\in \Omega_-\cap u_{\varepsilon}\cap \mathbb{C}^+$, there exists an integer $M>0$ and a constant $C > 0$, such that if $m\ge M$, there is a zero $z_k=z_k(\xi_k)$ of $\pi_m$ satisfying
$$|\,\xi_k^*-\xi_k\,|\le \frac{C}{m^{1/4}}.$$
\end{theorem}

The paper is organized as follows. Section \ref{sec-orthogonality} contains the orthogonality relation and the associated Riemann-Hilbert problem for the family of Jacobi Polynomials $\{\pi_m\}$ defined in \eqref{Eq def pi_m}. Section \ref{sec-RH Analysis} contains the asymptotic analysis of the Riemann Hilbert problem. The global asymptotics of the polynomials are determined in terms of a $2 \times 2$ matrix, whose $(1,1)$ entry provides the asymptotics of the polynomial $\pi_m$. Finally, Section \ref{sec zeros} contains the explicit asymptotic formulas in the regions where the zeros are located. Theorem \ref{Thm Asymptotics} is proved here, and then it is used to derive the asymptotic results for the location of the zeros, establishing Theorems \ref{Thm Zeros out} and  \ref{Thm Zeros in}.

\section{Orthogonality and Riemann Hilbert formulation}
\label{sec-orthogonality}
Following the approach of Kuijlaars, Martinez-Finkelshtein, and Orive \cite{kuijlaars2005orthogonality}, weight functions are  defined continuously over a closed contour \( \Gamma \) encircling the points \( +1 \) and \( -1 \). This yields an orthogonality relation along the contour for a family of Jacobi polynomials. An associated RHP then shows that this orthogonality relation characterizes the family of polynomials.
\subsection{Orthogonality in the complex plane}
\label{subsec complex orthogonality}
\begin{figure}[htbp]
    \centering
    \begin{tikzpicture}[scale=0.8]
        \draw[<->] (-3,0) -- (6,0) node[below]{};
        \node at (2,0) [circle,fill,inner sep=1pt,label=below:0] {};
        \node at (0,0) [circle,fill,inner sep=1pt,label=below:-1] {};
        \node at (4,0) [circle,fill,inner sep=1pt,label=below:1] {};
        \node at (2,1.7) [above] {$\Gamma$}; 

        \draw[line width=1pt,blue!50!, decoration={
            markings,
            mark=at position 0.2 with {\arrow[line width=1.2pt]{>}},
            mark=at position 0.6 with {\arrow[line width=1.2pt]{>}},
            mark=at position 0.8 with {\arrow[line width=1.2pt]{>}}
        }, postaction={decorate}] 
        
        (2,0) to[out=-60,in=-90] (5,0)
        to[out=90,in=90] (-1,0) to[out=-90,in=-120] (2,0); 
        
    \end{tikzpicture}
    \caption{Original contour $\Gamma$.}
    \label{fig_Gamma_original}
\end{figure}
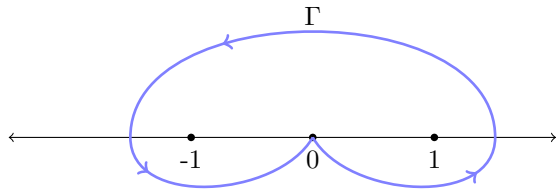

Let $\Gamma$ be a counterclockwise oriented contour encircling $[-1,1]$ as in Figure \ref{fig_Gamma_original}. Note that it is permitted to touch but not cross the interval $(-1,1]$.

For any non-negative integer $n$ consider the weight function 
$$w_n(t)=(1-t)^{n+\frac{1}{2}}(1+t)^{-n-\frac{1}{2}}.$$ 
This can be defined continuously on the contour $\Gamma$ considering the branch cut on $[-1,1]$, starting with the positive value $1$ at zero in the lower side of the branch cut, and then moving forward on $\Gamma$ in the positive orientation. More specifically, define the complex valued functions $(1+z)^{\frac{1}{2}}$ and $(1-z)^{\frac{1}{2}}$, with branch cuts opening in the the right direction, as follows:
\begin{equation}\label{Eq Def sqrt(1+-z)}
\begin{split}
(1+z)^{\frac{1}{2}}&=|(1+z)^{\frac{1}{2}}|e^{\frac{i\theta}{2}}, \,\,\, \,\,\,\,\,\, \,\,\,\,0\le \theta<2\pi,\\
 (1-z)^{\frac{1}{2}}& =|(1-z)^{\frac{1}{2}}|e^{\frac{i(\phi+\pi)}{2}}, \,\,\,\,0\le \phi<2\pi,
 \end{split} 
\end{equation}
with $\phi=\textnormal{Arg}(z-1)$ and $\theta=\textnormal{Arg}(z-(-1))$. Then
\begin{equation}\label{EQ def hat w}
     \widetilde{w}(z)=\frac{(1-z)^{\frac{1}{2}}}{(1+z)^{\frac{1}{2}}}
\end{equation}
is analytic in $\mathbb{C}\setminus [-1,1]$, with
$$\widetilde{w}_{-}(0):= \lim_{\substack{z \to 0 \\ \operatorname{Im}(z) < 0}} \widetilde{w}(z) = 1 > 0.$$
The weight function $w_n(t)$ is then defined continuously on $\Gamma$ via
\begin{equation}\label{Eqtn_wn(z)}
w_n(t)=\left(\frac{1-t}{1+t}\right)^n \widetilde{w}(t).
\end{equation}
For fixed $n\ge 0$,  \cite[Thm. 2.1]{ kuijlaars2005orthogonality} gives an orthogonality 
relation for the family of Jacobi polynomials 
$$\{p_k(z)\}_{k\ge 0}=\{P_k^{\left(n+\frac{1}{2},-n-\frac{1}{2}\right)}(z)\}_{k\ge 0},$$
over the contour $\Gamma$, with respect to the weight function $w_n(t)$; that is,
\begin{equation}\label{Eq. orthogonality}
    \int_{\Gamma}t^jp_k(t)w_n(t)\, dt=h(j,k) \delta_{jk},\hspace{0.8cm}\text{for}\,\,  j=0,1,2,\ldots, k,
\end{equation}
where
\begin{equation}\label{Eq values h(n,m)}
    h(j,k)=\frac{-\pi^2\, 2^{k+2}}{\Gamma(2k+2)\Gamma\left(-k-j-\frac{1}{2}\right)\Gamma\left(-k+j+\frac{1}{2}\right)}.
\end{equation}
Here $\Gamma(\,\cdot\,)$ is the classical Eulerian gamma function \cite{andrews1999special}.
\begin{remark}
    \begin{enumerate}
        \item The values of $h(j,k)$ in \eqref{Eq values h(n,m)} are half of those from Theorem 2.1 in \cite{kuijlaars2005orthogonality}. This is due to the fact that, in the current case, a single loop around the points $\pm 1$ suffices for the orthogonality relation to hold.
       \item The family of polynomials $P_m(z)$ defined in \eqref{jaco-explicit2} corresponds to the diagonal of the two-dimensional family $\{P_k^{(n+1/2,-n-1/2)}\}_{k,n\ge 0}.$ That is   
       $$P_m(z)=P_m^{\left(m+\frac{1}{2},-m-\frac{1}{2}\right)}(z).$$
    \end{enumerate}
\end{remark}
\subsection{Riemann Hilbert Problem for Orthogonal Polynomials}
\label{subsec RHP for OP}
The RHP for orthogonal polynomials was first introduced  by Fokas, Its and Kitaev \cite{fokas1992isomonodromy} in the context of orthogonality on the real axis. In the case of Jacobi polynomials with non classical parameters, this was extended using orthogonality on a complex contour by Kuijlaars, Mart\'inez-Finkelshtein, Mart\'inez-Gonz\'alez and Orive in \cite{kuijlaars2005orthogonality}. The Riemann Hilbert formulation of Jacobi polynomials was then used to describe their asymptotic behavior and zeros, see \cite{ kuijlaars2005orthogonality,  kuijlaars2004strong, martinez2005riemann}. Previous work, using potential theory, appears in \cite{martinez2000zeros}. 

For each $n\ge 0$ we now formulate the RHP corresponding to the family of polynomials $\{P_k^{(n+1/2,-n-1/2)}(z)\}_{k \ge 0}$.

\begin{figure}[H]
    \centering
    \includegraphics[width=0.4\textwidth]{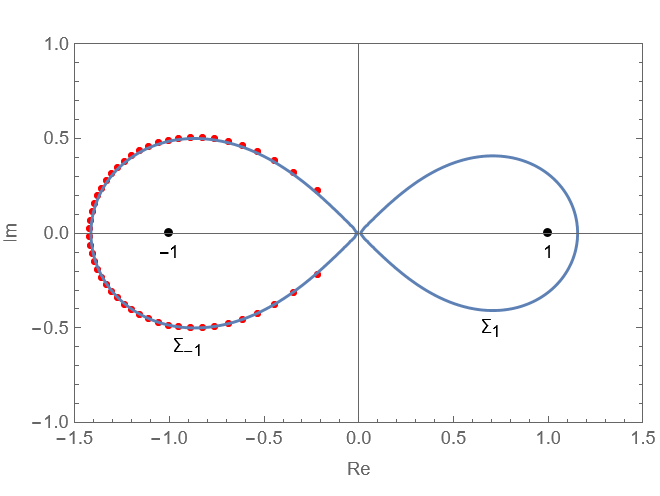} 
    \caption{Contour $\Sigma$ and zeros of $P_{50}(z)=P_{50}^{\left(50+\frac{1}{2},-50-\frac{1}{2}\right)}(z)$.}
    \label{fig: Zeros_lemniscate_p50}
\end{figure}

\begin{RHP problem}\label{RHP Y}
Let $\Gamma$ be the contour in Figure \ref{fig_Gamma_original}, and let $n\in \mathbb{Z}$, $n\ge 0$. The problem is to determine a $2\times 2$ matrix-valued function $Y(z): \mathbb{C}\setminus \Gamma \rightarrow\mathbb{C}^{2\times 2}$ satisfying the following conditions:
\begin{equation}\label{Eq. RHP Y}
\begin{cases}
Y(z) \,\,\,\text{is analytic in}\,\, \mathbb{C}\setminus \Gamma,\\
Y_+(t)=Y_-(t)
     \begin{pmatrix}
       1 & w_n(t)\\
       0 & 1 \\
    \end{pmatrix},\hspace{1cm}t\in \Gamma,\\
Y(z)=\left(I+O\left(z^{-1}\right)\right)
    \begin{pmatrix}
                    z^k & 1 \\
                    0 & z^{-k}\\
    \end{pmatrix},
    \hspace{1cm}
    z\rightarrow\infty.
\end{cases}    
\end{equation}
\end{RHP problem}
The RHP \ref{RHP Y} has a unique solution, expressed in terms of Jacobi polynomials, see Section 3 in \cite{kuijlaars2005orthogonality}. This is given next.
\begin{proposition}\label{RHP Y sol polynomials}
    The unique solution of the RHP \eqref{RHP Y} is given by
\begin{center}
  Y(z)\,=$\begin{pmatrix}
       \pi_k(z) & \frac{1}{2\pi i} \int_{\Gamma}\frac{\pi_k(t)w_n(t)}{t-z}\, dt\\
       c_{k-1} \pi_{k-1}(z) & \frac{c_{k-1}}{2\pi i} \int_{\Gamma}\frac{\pi_{k-1}(t)w_n(t)}{t-z}\, dt\\
    \end{pmatrix}.$
\end{center}
Here $\pi_k(z)$ denotes the monic Jacobi polynomial of degree $k$, with respect to the weight $w_n(t)$ defined in equation \eqref{Eqtn_wn(z)}. The constant in the second row is \linebreak$c_k=\frac{-2\pi i}{h(n,k)}$, with $h$ given in \eqref{Eq values h(n,m)}.
\end{proposition}
The contour $ \Gamma $ can be deformed as long as it crosses neither the branch cut on $ [-1,1] $ nor the singular point $ z = -1 $. The contour $\Gamma$ will be deformed into the contour $\Sigma$ defined in \eqref{Eq Contour Sigma}. As stated in Corollary \ref{Corollary Asymp Zeros}, it will be shown that the zeros of the polynomials $P_m(z)=P_m^{(m+1/2,-m-1/2)}(z)$ accumulate on $\Sigma_{-1}$ as $m\rightarrow\infty$, see Figure \ref{fig: Zeros_lemniscate_p50} for an illustration.
\section{Riemann Hilbert Analysis}
\label{sec-RH Analysis}
In this section, the asymptotic solution of RHP \ref{RHP Y} is established in the special case $n=m$. The process involves applying a sequence of transformations to the original problem for $Y(z)$ (see \eqref{Eq. RHP Y}), turning it into a problem for $U(z)$ (see \eqref{RHP U}):
$$ Y \longrightarrow T \longrightarrow U, $$
along with contour deformations whose role is to simplify the problem. It turns out that the origin $z=0$ requires a special analysis.

This section is organized as follows. Subsection \ref{sec Outside Solution} introduces the global transformations leading to the RHP for $U(z)$. Then an asymptotic solution, away from the origin, to the problem for $U(z)$ is established. Subsection \ref{sec-RHP PC} introduces the Parabolic Cylinder RHP, which is later used to construct the local solution near the origin. Subsection \ref{sec local parametrix} describes the local RHP \ref{RHP local} and establishes a local solution to the problem for $U(z)$, in terms of the solution to the Parabolic Cylinder RHP. Finally, in Subsection \ref{sec Global solution} we use the local and outer solutions to establish global uniform asymptotics for $Y(z)$ in the entire complex plane.

\subsection{Solution away from the origin}
\label{sec Outside Solution}
Throughout the paper, $\sigma_3$ is the matrix 
\begin{equation*}
   \sigma_3= \begin{pmatrix}
        1&0\\
        0&-1
    \end{pmatrix}.
\end{equation*}
For any complex number $z$, the expression $z^{\sigma_3}$ is given by  
\begin{equation*}
   z^{\sigma_3} =\begin{pmatrix}
        z&0\\
        0&z^{-1}
    \end{pmatrix}.
\end{equation*}

First introduce the transformation $T(z)=Y(z)H(z)$, where $Y(z)$ is the solution to RHP \ref{RHP Y}, and
\begin{equation*}
    H(z)= \begin{cases} 
      (1-z)^{m\sigma_3}   & \text{if $z$ is inside $\Sigma_{-1}$}, \\
      (1+z)^{-m\sigma_3} &\text{otherwise.}
   \end{cases}
\end{equation*}
RHP \ref{RHP Y} now yields the RHP for $T(z)$:
\begin{equation}\label{RHP T}
\begin{cases}
T(z) \,\,\,\text{is analytic in}\,\, \mathbb{C}\setminus \Sigma,\\
T_+(t)=T_-(t)J_T(t)\,\,\text{if}\,\, t\in \Sigma,\\
T(z)=I+O\left(\frac{1}{z}\right)\,\,\text{as}\,\, z\rightarrow\infty, 
\end{cases}    
\end{equation}
with
\begin{equation*}
    J_T(t)= \begin{cases} 
         \begin{pmatrix}
       (1-t^2)^m & \widetilde{w}(t)\\
       0 & (1-t^2)^{-m} \\
    \end{pmatrix},& \text{if $t\in \Sigma_{-1}$,}\\[0.4cm]
         \begin{pmatrix}
       1 & (1-t^2)^m\, \widetilde{w}(t)\\
       0 & 1 \\
    \end{pmatrix},&\text{if $t\in \Sigma_{1}$,}
   \end{cases}
\end{equation*}
and the function $\widetilde{w}(z)$, defined in \eqref{EQ def hat w}, is analytic in $\mathbb{C}\setminus [-1,1]$.

Let 
\begin{equation*}
    L(t)=\begin{pmatrix}
       1& 0\\
      \frac{1}{(1-t^2)^{m}\,\widetilde{w}(t)} & 1
    \end{pmatrix},\hspace{0.5cm}
    M(t)=\begin{pmatrix}
       0& \widetilde{w}(t)\\
       -\frac{1}{\widetilde{w}(t)}& 0
       \end{pmatrix},
\end{equation*}
\begin{equation*}
       R(t)=\begin{pmatrix}
       1& 0\\
      \frac{(1-t^2)^{m}}{\widetilde{w}(t)}& 1
    \end{pmatrix}\,\,\, \text{and}\,\,\, A(t)=\begin{pmatrix}
       1 & (1-t^2)^m\, \widetilde{w}(t)\\
       0 & 1 \\
    \end{pmatrix}.
\end{equation*}
A direct calculation gives the factorization
\begin{equation*}
\begin{pmatrix}
       (1-t^2)^m& \widetilde{w}(t)\\
       0 & (1+t^2)^{-m}
    \end{pmatrix}= L(t)M(t)R(t).
\end{equation*}

This implies that the transformation $U(z)=T(z)K(z)$, with
\begin{equation*}
    K(z)= \begin{cases} 
        \hspace{2cm}I & \text{if $z\in \Omega_0 \cup \Omega_{\infty}$}, \\
         \begin{pmatrix}
       1& 0\\
      \frac{(1-z^2)^{m}}{ \widetilde{w}(z)} & 1
    \end{pmatrix}^{-1}&\text{if $z\in \Omega_+$},\\
        \begin{pmatrix}
       1& 0\\
      \frac{1}{(1-z^2)^{m}\,\widetilde{w}(z)}  & 1
    \end{pmatrix} &\text{if $z\in \Omega_-$},
   \end{cases}
\end{equation*}
yields the following RHP for \(U(z)\):
\begin{equation}\label{RHP U}
\begin{cases}
U(z) \,\,\,\text{is analytic in}\,\, \mathbb{C}\setminus \Sigma_U=\Sigma\cup\Sigma_r\cup\Sigma_{1/r},\\
 U_+(t)=U_-(t)J_U(t)\hspace{0,5cm} \text{if}\,\, t\in \Sigma_U,\\
U(z)=I+O\left(\frac{1}{z}\right)\hspace{0,5cm}\text{as}\,\, z\rightarrow\infty, 
\end{cases}    
\end{equation}
with
\begin{equation}\label{Eq jump U}
    J_U(z)= \begin{cases} 
      M(t)   & \text{if $t\in \Sigma_{-1}$}, \\
      R(t)   & \text{if $t\in \Sigma_{r}$}, \\
      L(t)   & \text{if $t\in \Sigma_{1/r}$}, \\
      A(t)   & \text{if $t\in \Sigma_{1}$}.
   \end{cases}
\end{equation}
For $z\in \Sigma_{r},\,\Sigma_{1/r}$, or $\Sigma_{1}$, and $z$ bounded away from $z=0$, the jump matrices $R$, $L$, and $A$ are close to the identity (for  $m$ is sufficiently large). This intuitively tells us that on these portions of the contour $\Sigma$, the jump relations for $U$ are almost negligible. Hence the leading order asymptotics of the solution to the RHP for $U(z)$ \eqref{RHP U}, are determined by the solution of the following RHP.

\begin{RHP problem}\label{RHP Outer}
Find a  piecewise analytic $2 \times 2$ matrix-valued function $N(z): \mathbb{C} \setminus \Sigma \rightarrow \mathbb{C}^{2\times 2}$ satisfying:
\begin{equation}\label{RHP N}
\begin{cases}
N(z)\,\, \text{is analytic in}\,\, \mathbb{C}\setminus \Sigma_{-1},\\
 N_+(t)=N_-(t)\,M(t), \,\,\text{if}\,\,\, t\in \Sigma_{-1},\\
N(z)=I+\mathcal{O}\left(\frac{1}{z}\right),\,\,\text{as}\,\, z\rightarrow\infty. 
\end{cases}    
\end{equation}

\end{RHP problem}
The solution to this problem is given by
\begin{equation}\label{Eq N explicit}
    N(z)=\begin{cases} 
        \begin{pmatrix}
        0&\frac{(1-z)^{\frac{1}{2}}}{z^{\frac{1}{2}}}\\
        -\frac{z^{\frac{1}{2}}}{(1-z)^{\frac{1}{2}}}& 0
    \end{pmatrix}  & \text{if $z$ is inside $\Sigma_{-1}$},\\
       \begin{pmatrix}
        \frac{(1+z)^{\frac{1}{2}}}{z^{\frac{1}{2}}}& 0 \\
        0 & \frac{z^{\frac{1}{2}}}{(1+z)^{\frac{1}{2}}}
    \end{pmatrix}  &\text{otherwise,}
   \end{cases}
\end{equation}
here $z^{\frac{1}{2}}$ is defined by $z^{\frac{1}{2}} :=\abs{z^{\frac{1}{2}}}e^{\frac{i\,\textnormal{Arg}(z)}{2}}, \,\,\,\,0\le\textnormal{Arg}(z)<2\pi$.

For $z$ bounded away from $z=0$, the asymptotics of the solution to the RHP in \eqref{RHP U}, as we will show, are determined by $N(z)$. The next step is a local analysis of the problem near the origin. First we introduce a RHP that will be used to establish the local solution of RHP in \eqref{RHP U}.
\subsection{The parabolic cylinder RHP}
\label{sec-RHP PC}
This subsection introduces a RHP whose solution will be used to construct a local parametrix near the origin. The problem is named after its explicit solution, which is expressed in terms of parabolic cylinder functions. The problem, as presented in \cite{bothner2012nonlinear}, is formulated as follows.

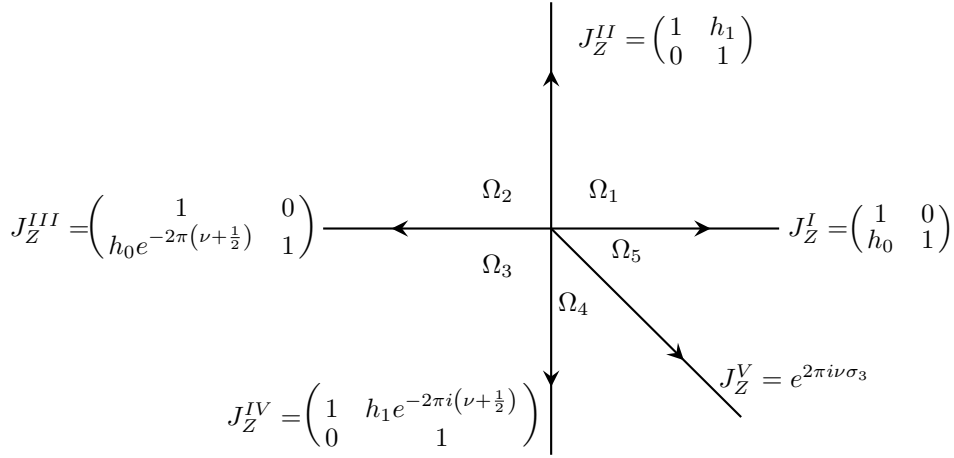
\begin{figure}[htbp]
    \centering
    \begin{tikzpicture}[
        myarrow/.style={
            thick,
            postaction={decorate,
            decoration={markings,mark=at position #1 with {\arrow[scale=1.5]{stealth}}}}
        },
        every matrix/.style={scale=0.5} 
    ]
        \draw[myarrow=0.7] (0,0) -- (0,3);
        \draw[myarrow=0.7] (0,0) -- (0,-3);
        \draw[myarrow=0.7] (0,0) -- (-3,0);
        \draw[myarrow=0.7] (0,0) -- (3,0);
        \draw[myarrow=0.7] (0,0) -- (2.5,-2.5);

        \node (eq) at (0.7,0.5) {$\Omega_1$};
        \node (eq) at (-0.7,0.5) {$\Omega_2$};
        \node (eq) at (0.3,-1) {$\Omega_4$};
        \node (eq) at (-0.7,-0.5) {$\Omega_3$};
        \node (eq) at (1,-0.3) {$\Omega_5$};

        \node (eq) at (3.2,-2) {$J_Z^{V}=e^{2\pi i\nu \sigma_3}$};
        
        \matrix (m2) at (4.6,0) [matrix of math nodes, left delimiter=(, right delimiter=), inner sep=0pt, column sep=10pt, row sep=3pt] {
            1 & 0  \\
           h_0& 1\\
        }; \node[anchor=east] at (m2.west) {$J_Z^{I}=\,\,$};
        
        \matrix (m3) at (2,2.5) [matrix of math nodes, left delimiter=(, right delimiter=), inner sep=0pt, column sep=10pt, row sep=3pt] {
            1 & h_1\\
           0 & 1 \\
        };
        \node[anchor=east] at (m3.west) {$J_Z^{II}=\,\,\,$};
        
        \matrix (m4) at (-4.6,0) [matrix of math nodes, left delimiter=(, right delimiter=), inner sep=0pt, column sep=10pt, row sep=3pt] {
            1 & 0\\
           h_0 e^{-2\pi\left(\nu+\frac{1}{2}\right)} & 1 \\
        };
        \node[anchor=east] at (m4.west) {$J_Z^{III}=\,\,$};

        \matrix (m3) at (-1.7,-2.5) [matrix of math nodes, left delimiter=(, right delimiter=), inner sep=0pt, column sep=10pt, row sep=3pt] {
            1 & h_1 e^{-2\pi i \left(\nu+\frac{1}{2}\right)}\\
           0 & 1 \\
        };
        \node[anchor=east] at (m3.west) {$J_Z^{IV}=\,\,\,$};
    \end{tikzpicture}
    \caption{Contour $\Sigma_Z$ and  jump matrices $J_Z(\eta)$.}
    \label{Contour P}
\end{figure}

\begin{RHP problem} \label{RHP parabolic}
  Let $\Sigma_Z$ be the contour consisting of five rays emanating from the origin, one at each of the angles $0$,  $\frac{\pi}{2}$, $\pi$, $\frac{3\pi}{4}$ and $-\frac{\pi}{4}$. Given $r_0\in \mathbb{C}$ with $|r_0|> 1$, set $\nu=- \frac{1}{2\pi}\log(1-|r_0|^2)$. Find a piecewise analytic $2\times 2$ matrix-valued function $Z(\eta): \mathbb{C}\setminus \Sigma_Z \rightarrow\mathbb{C}^{2\times 2}$ satisfying the following:
\begin{enumerate}
    \item $Z(\eta)$  is analytic in $\mathbb{C}\setminus \Sigma_Z$,
    \item $Z_+(\eta)=Z_-(\eta)J_Z(\eta)$, for $\eta\in \Sigma_Z,$
    \item As $\eta\rightarrow\infty$,    
         \begin{equation}\label{Eq Asymptotics PC}
          Z(\eta)=\frac{\eta^{-\frac{\sigma_3}{2}}}{\sqrt{2}}\left[ \begin{pmatrix}
            1+\frac{(\nu+1)(\nu+2)}{2\eta^2}&1-\frac{\nu(\nu-1)}{2\eta^2}\\
            1+\frac{(\nu+1)(\nu-2)}{2\eta^2}&-1+\frac{\nu(\nu+3)}{2\eta^2}
        \end{pmatrix}+\mathcal{O}\left(\frac{1}{\eta^4}\right)\right] e^{\left(\frac{\eta^2}{4}-\left(\nu+\frac{1}{2}\right)\ln \eta\right)\sigma_3},
      \end{equation}
\end{enumerate}
where $J_Z(\eta)$ is given in Figure \ref{Contour P} and 
\begin{equation*}
    h_0=-i\frac{\sqrt{2\pi}}{\Gamma(\nu+1)}, \hspace{0.5cm} h_1=\frac{\sqrt{2\pi}}{\Gamma(-\nu)}e^{\pi i \nu}.
\end{equation*}  
\end{RHP problem}

The solution to RHP \ref{RHP parabolic} is described below (see \cite{bothner2012nonlinear} for more detail). It is given in terms of the parabolic cylinder function $D_\nu(\eta)$. This is the unique solution of the equation
\begin{equation*}
    \frac{d^2}{d\eta^2} D_{\nu}(\eta)+ \left(\frac{1}{2}-\frac{\eta^2}{2}+a\right)D_{\nu}(\eta)=0,
\end{equation*}
with the following asymptotic behavior as $\eta \rightarrow \infty$ (see \cite{Bateman1953}, Section 8.4),
\begin{equation*}\label{asymp parabolic cylinder}
    D_\nu(\xi)=\begin{cases}
    e^{-\frac{\eta^2}{4}}\eta^{\nu}\left(1-\frac{\nu(\nu-1)}{2\eta^2}+\mathcal{O}\left(\frac{1}{\eta^4}\right)\right),&  -\frac{3\pi}{4} <\textnormal{Arg}(\eta) <  \frac{3\pi}{4};\\
     e^{-\frac{\eta^2}{4}}\eta^{\nu}\left(1-\frac{\nu(\nu-1)}{2\eta^2}+\mathcal{O}\left(\frac{1}{\eta^4}\right)\right) & \\
     \hspace{1cm}-\frac{\sqrt{2\pi}}{\Gamma(-\nu)}e^{\nu i}e^{\frac{\eta^2}{4}}\eta^{-\nu-1}\left(1+\frac{(\nu+1)(\nu+2)}{2\eta^2}+\mathcal{O}\left(\frac{1}{\eta^4}\right)\right),&  \frac{\pi}{4} <\textnormal{Arg}(\eta) <  \frac{5\pi}{4};\\
     e^{-\frac{\eta^2}{4}}\eta^{\nu}\left(1-\frac{\nu(\nu-1)}{2\eta^2}+\mathcal{O}\left(\frac{1}{\eta^4}\right)\right) & \\
     \hspace{1cm}-\frac{\sqrt{2\pi}}{\Gamma(-\nu)}e^{-\nu i}e^{\frac{\eta^2}{4}}\eta^{-\nu-1}\left(1+\frac{(\nu+1)(\nu+2)}{2\eta^2}+\mathcal{O}\left(\frac{1}{\eta^4}\right)\right),&  -\frac{5\pi}{4} <\textnormal{Arg}(\eta) <  -\frac{\pi}{4}.
    \end{cases}
\end{equation*}
For $\eta\in \mathbb{C}$, define
\begin{equation*}
    Z_0(\eta) = 2^{-\sigma_3}
    \begin{pmatrix}
        D_{-\nu-1}(i\eta) & D_{\nu}(\eta) \\
        \frac{d}{d\eta} D_{-\nu-1}(i\eta) & \frac{d}{d\eta} D_{\nu}(\eta)
    \end{pmatrix}
    \begin{pmatrix}
        e^{i\frac{\pi}{2}(\nu+1)} &0\\
        0&1
    \end{pmatrix},
\end{equation*}
and define recursively 
\begin{equation}\label{Eq. def Z_n}
    Z_{n+1}(\eta)=Z_n(\eta)H_n, \hspace{0.5cm}\text{$n=0,1,2,3$,}
\end{equation}
where
\begin{equation*}
    H_0=\begin{pmatrix}
        1&0\\
        h_0& 1
    \end{pmatrix}, \hspace{0.3cm} 
 H_1=\begin{pmatrix}
     1&h_1\\
        0& 1
 \end{pmatrix}, \hspace{0.3cm}h_0=-i\frac{\sqrt{2\pi}}{\Gamma(\nu+1)},\hspace{0.3cm}h_1=-i\frac{\sqrt{2\pi}}{\Gamma(-\nu)}e^{i\pi\nu},
\end{equation*}
and for $n=0,1$,
\begin{equation*}
 H_{n+2}=e^{i\pi\left(\nu+\frac{1}{2}\right)\sigma_3}\,H_{n}\, e^{-i\pi\left(\nu+\frac{1}{2}\right)\sigma_3}.   
\end{equation*}
The solution of RHP \ref{RHP parabolic} is given by, 
\begin{equation}\label{Eq sol RHP for P}
    Z(\eta) = 
    \begin{cases}
        Z_0(\eta), & -\frac{\pi}{4} <\textnormal{Arg}(\eta) < 0;\\        
        Z_1(\eta), &  0 < \textnormal{Arg}(\eta) < \frac{\pi}{2};\\
        Z_2(\eta), &  \frac{\pi}{2} < \textnormal{Arg}(\eta) < \pi;\\
        Z_3(\eta), &  \pi < \textnormal{Arg}(\eta) < \frac{3\pi}{2};\\
        Z_4(\eta), &  \frac{3\pi}{2} < \textnormal{Arg}(\eta) < \frac{7\pi}{4}.
    \end{cases}
\end{equation}
Details may found in \cite{bothner2012nonlinear}.
\subsection{Local solution to RHP 3.2 near the origin}
\label{sec local parametrix}

The parabolic cylinder RHP is now used to construct a local parametrix $X(z)$ near the origin. This is a locally defined solution that approximates the true solution $U(z)$ in a  small disk centered at $z=0$, say $u_{\varepsilon}= \{z\in \mathbb{C}: |z|<\varepsilon\}$, where the matrix $N(z)$ in \eqref{Eq N explicit} is not a good approximation of $U(z)$. The local parametrix must be constructed so that:
\begin{itemize}
    \item $X(z)$ satisfies the jumps relations of $U(z)$ exactly, inside  $u_{\varepsilon}$, see \eqref{Eq jump U},
    \item $X(z)$ matches $N(z)$ on the boundary of $u_{\varepsilon}$ up to order $1/m$.
\end{itemize}
More precisely, let $\Sigma_{\varepsilon}$ denote the boundary of the disk $u_{\varepsilon}$, then consider:
\begin{RHP problem}\label{RHP local}
Let $\varepsilon>0$ be small enough, find a  piecewise analytic $2\times 2$  matrix-valued function $X(z): \bar{u}_{\varepsilon}\setminus \Sigma_U \rightarrow\mathbb{C}^{2\times 2}$ satisfying:
\begin{equation*}
\begin{cases}
X(z) \text{ is analytic in $u_{\varepsilon}\setminus \Sigma_U$ and continuous on $\bar{u}_{\varepsilon}\setminus \Sigma_U$}, \\
X_+(t) = X_-(t) J_U(t), \text{ for } t \in u_{\varepsilon}\cap \Sigma_U, \\
\text{There is a constant $C > 0$ such that for every $z \in\Sigma_{\varepsilon}\setminus \Sigma_U$, we have}\\
\hspace{3cm}\|N(z)X^{-1}(z) - I\| \le \frac{C}{m},\hspace{0.3cm}
\text{for any matrix norm $||\,\cdot\,||$.}
\end{cases}    
\end{equation*}

\end{RHP problem}
We first seek a solution that satisfies the jump relations, temporarily disregarding the boundary condition of the problem on the disk $u_{\varepsilon}$. A sequence of transformations: 
\begin{equation*}
X=X^{(1)}\rightarrow X^{(2)}\rightarrow X^{(3)}\rightarrow X^{(4)} \rightarrow X^{(5)},
\end{equation*}
turns the problem for $X=X^{(1)}$ into a problem for $X^{(5)}$, having the same jump relations as the Parabolic Cylinder RHP \ref{RHP parabolic}.  In subsection \ref{subsec Adjust local sol} the asymptotics of the Parabolic Cylinder RHP (see \eqref{Eq Asymptotics PC}) are used to obtain a solution satisfying the desired boundary conditions.

First consider the function
\begin{equation*}
    \widehat{w}(z)= \begin{cases}
                \widetilde{w}(z) & \text{if $z\in \mathbb{C}_+$},\\
                -\widetilde{w}(z) & \text{if $z\in \mathbb{C}_-$},
             \end{cases}
\end{equation*}
which is analytic in $\mathbb{C}\setminus\{(-\infty,-1]\cup[1,\infty)\}$, since the boundary values $\widehat{w}_+(t)$  and $\widehat{w}_-(t)$ agree on this interval. The disk $u_{\varepsilon}$ can be chosen so that $\widehat{w}(z)$ is analytic inside it. A further condition on $u_{\varepsilon}$ will appear later, see \eqref{Eq def map lambda}. The analysis below is all performed in such a disk.

\begin{figure}[H]
    \centering
    \begin{tikzpicture}[
        myarrow/.style={
            thick,
            postaction={decorate,
            decoration={markings,mark=at position #1 with {\arrow[scale=1.5]{stealth}}}}
        },
        every matrix/.style={scale=0.4} 
    ]
        \draw[myarrow=0.5] (0,0) -- (3,-2);
        \draw[myarrow=0.5] (3,2) -- (0,0);
        
        \draw[myarrow=0.5] (0,0) -- (0,4);
        \draw[myarrow=0.5] (0,0) -- (-3,3);
        \draw[myarrow=0.5] (0,0) -- (-4,2);

        \draw[myarrow=0.5] (-4,-2) -- (0,0) ;
        \draw[myarrow=0.5] (-3,-3) -- (0,0) ;
        \draw[myarrow=0.5] (0,-4) -- (0,0) ;
        \draw[dashed] (-6,0) -- (6,0);
        
        \matrix (m1) at (4,1.4) [matrix of math nodes, left delimiter=(, right delimiter=), inner sep=0pt, column sep=10pt, row sep=3pt] {
            1 & (1-t^2)^m\widehat{w} \\
           0 & 1 \\
        };
        
        \matrix (m2) at (4,-1.4) [matrix of math nodes, left delimiter=(, right delimiter=), inner sep=0pt, column sep=10pt, row sep=3pt] {
            1 & -(1-t^2)^m\widehat{w} \\
           0& 1\\
        };
        \matrix (m3) at (-4,3) [matrix of math nodes, left delimiter=(, right delimiter=), inner sep=0pt, column sep=10pt, row sep=3pt] {
            0 & \widehat{w}\\
           -\widehat{w}^{-1} & 0 \\
        };
        \matrix (m4) at (1.15,3.4) [matrix of math nodes, left delimiter=(, right delimiter=), inner sep=0pt, column sep=10pt, row sep=3pt] {
            1 & 0\\
            \frac{\widehat{w}^{-1}}{(1-t^2)^m} & 1 \\
        };
        \matrix (m4) at (-5,1.4) [matrix of math nodes, left delimiter=(, right delimiter=), inner sep=0pt, column sep=10pt, row sep=3pt] {
            1 & 0\\
            (1-t^2)^m\widehat{w}^{-1} & 1 \\
        };
        \matrix (m4) at (-5,-1.4) [matrix of math nodes, left delimiter=(, right delimiter=), inner sep=0pt, column sep=10pt, row sep=3pt] {
            1 & 0\\
            -(1-t^2)^m\widehat{w}^{-1} & 1 \\
        };
        \matrix (m4) at (1.15,-3.4) [matrix of math nodes, left delimiter=(, right delimiter=), inner sep=0pt, column sep=10pt, row sep=3pt] {
            1 & 0\\
            \frac{-\widehat{w}^{-1}}{(1-t^2)^m} & 1 \\
        };
        \matrix (m3) at (-4,-3) [matrix of math nodes, left delimiter=(, right delimiter=), inner sep=0pt, column sep=10pt, row sep=3pt] {
            0 & -\widehat{w}\\
           \widehat{w}^{-1} & 0 \\}; 
           
        \node (eq) at (2,0.5) {$\widetilde{\Omega}_{I.1}$};
        \node (eq) at (2,-0.5) {$\widetilde{\Omega}_{I.2}$};
        \node (eq) at (0.7,2) {$\widetilde{\Omega}_{II}$};
        \node (eq) at (0.7,-2) {$\widetilde{\Omega}_{VIII}$};
        \node (eq) at (-1,2) {$\widetilde{\Omega}_{III}$};
        \node (eq) at (-1,-2) {$\widetilde{\Omega}_{VII}$};
        \node (eq) at (-2.5,1.8) {$\widetilde{\Omega}_{IV}$};
        \node (eq) at (-2.5,-1.8) {$\widetilde{\Omega}_{VI}$};
        \node (eq) at (-3,0.5) {$\widetilde{\Omega}_{V.1}$};
        \node (eq) at (-3,-0.5) {$\widetilde{\Omega}_{V.2}$}; 

        \node (eq) at (3.2,-2) {$\Sigma_1$};
        \node (eq) at (3.2,2) {$\Sigma_1$};
        \node (eq) at (-2.5,3) {$\Sigma_{-1}$};
        \node (eq) at (-2.4,-3) {$\Sigma_{-1}$};
        \node (eq) at (-0.4,-3.8) {$\Sigma_{1/r}$};
        \node (eq) at (-0.4,3.8) {$\Sigma_{1/r}$};
        \node (eq) at (-4.2,-2.1) {$\Sigma_{r}$};
         \node (eq) at (-4.2,2.1) {$\Sigma_{r}$};
    \end{tikzpicture}
    \caption{Jump matrices for $X^{(1)}(z)$.}
    \label{Contour X1}
\end{figure}
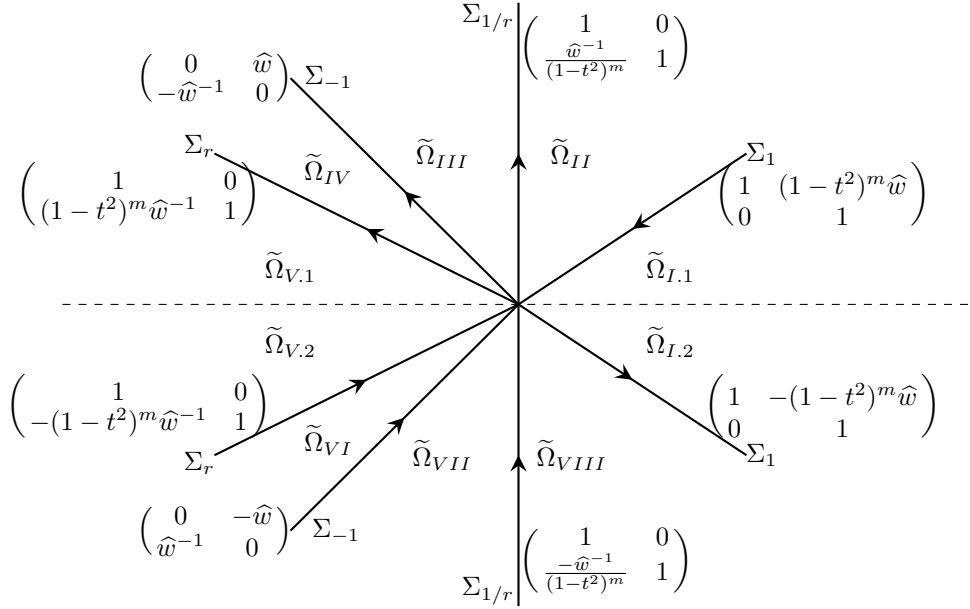

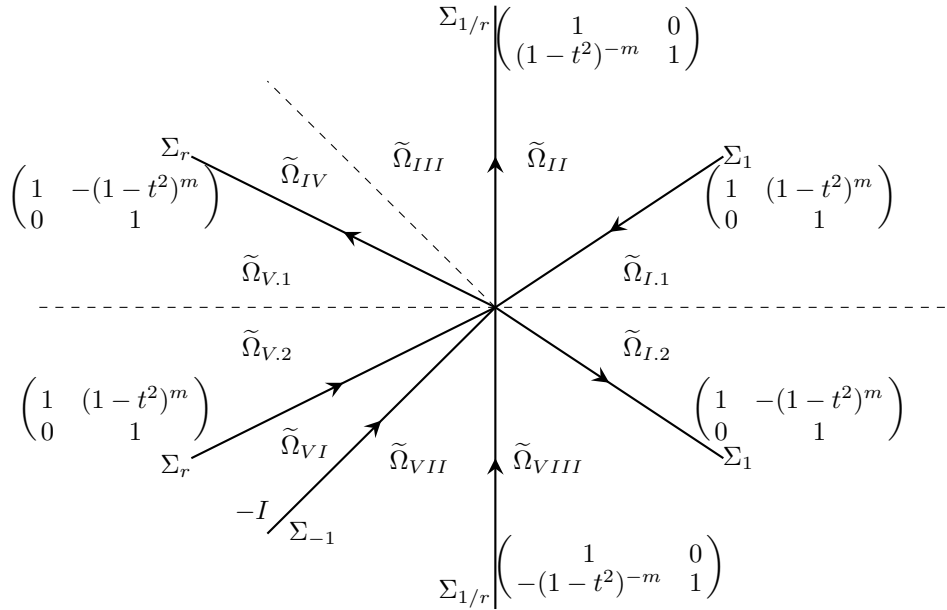
\begin{figure}[H]
    \centering
    \begin{tikzpicture}[
        myarrow/.style={
            thick,
            postaction={decorate,
            decoration={markings,mark=at position #1 with {\arrow[scale=1.5]{stealth}}}}
        },
        every matrix/.style={scale=0.4} 
    ]
         \draw[myarrow=0.5] (0,0) -- (3,-2);
        \draw[myarrow=0.5] (3,2) -- (0,0);
        
        \draw[myarrow=0.5] (0,0) -- (0,4);
        \draw[dashed] (0,0) -- (-3,3);
        \draw[myarrow=0.5] (0,0) -- (-4,2);

        \draw[myarrow=0.5] (-4,-2) -- (0,0) ;
        \draw[myarrow=0.5] (-3,-3) -- (0,0) ;
        \draw[myarrow=0.5] (0,-4) -- (0,0) ;
        \draw[dashed] (-6,0) -- (6,0);
        
        \matrix (m1) at (4,1.4) [matrix of math nodes, left delimiter=(, right delimiter=), inner sep=0pt, column sep=10pt, row sep=3pt] {
            1 & (1-t^2)^m \\
           0 & 1 \\
        };
        
        \matrix (m2) at (4,-1.4) [matrix of math nodes, left delimiter=(, right delimiter=), inner sep=0pt, column sep=10pt, row sep=3pt] {
            1 & -(1-t^2)^m \\
           0& 1\\
        };
        \matrix (m4) at (1.35,3.5) [matrix of math nodes, left delimiter=(, right delimiter=), inner sep=0pt, column sep=10pt, row sep=3pt] {
            1 & 0\\
            (1-t^2)^{-m}& 1 \\
        };
        \matrix (m4) at (-5,1.4) [matrix of math nodes, left delimiter=(, right delimiter=), inner sep=0pt, column sep=10pt, row sep=3pt] {
            1 & -(1-t^2)^m\\
            0& 1 \\
        };
        \matrix (m4) at (-5,-1.4) [matrix of math nodes, left delimiter=(, right delimiter=), inner sep=0pt, column sep=10pt, row sep=3pt] {
            1 & (1-t^2)^m\\
            0 & 1 \\
        };
        \matrix (m4) at (1.5,-3.5) [matrix of math nodes, left delimiter=(, right delimiter=), inner sep=0pt, column sep=10pt, row sep=3pt] {
            1 & 0\\
            -(1-t^2)^{-m} & 1 \\
        };
        \node (eq) at (-3.2,-2.7) {$-I$}; 
           
        \node (eq) at (2,0.5) {$\widetilde{\Omega}_{I.1}$};
        \node (eq) at (2,-0.5) {$\widetilde{\Omega}_{I.2}$};
        \node (eq) at (0.7,2) {$\widetilde{\Omega}_{II}$};
        \node (eq) at (0.7,-2) {$\widetilde{\Omega}_{VIII}$};
        \node (eq) at (-1,2) {$\widetilde{\Omega}_{III}$};
        \node (eq) at (-1,-2) {$\widetilde{\Omega}_{VII}$};
        \node (eq) at (-2.5,1.8) {$\widetilde{\Omega}_{IV}$};
        \node (eq) at (-2.5,-1.8) {$\widetilde{\Omega}_{VI}$};
        \node (eq) at (-3,0.5) {$\widetilde{\Omega}_{V.1}$};
        \node (eq) at (-3,-0.5) {$\widetilde{\Omega}_{V.2}$}; 

        \node (eq) at (3.2,-2) {$\Sigma_1$};
        \node (eq) at (3.2,2) {$\Sigma_1$};
        \node (eq) at (-2.4,-3) {$\Sigma_{-1}$};
        \node (eq) at (-0.4,-3.8) {$\Sigma_{1/r}$};
        \node (eq) at (-0.4,3.8) {$\Sigma_{1/r}$};
        \node (eq) at (-4.2,-2.1) {$\Sigma_{r}$};
         \node (eq) at (-4.2,2.1) {$\Sigma_{r}$};
    \end{tikzpicture}
    \caption{Jump matrices for $X^{(2)}(z)$.}
    \label{Contour X2}
\end{figure}

The jumps of the local RHP \ref{RHP local} can be expressed in terms of the function $\widehat{w}(z)$. Denote the solution to this problem as $X^{(1)}(z)$, with its corresponding jump matrices represented by $J_{X^{(1)}}(t)$. The labels of the regions in the new partition appears in the Figure \ref{Contour X1}.
\noindent
\textbf{First transformation:}  To remove the term $\widehat{w}(z)$ in the jump matrices and collapse the jump on the upper ray of $\Sigma_{-1}\cap u_{\varepsilon}$ to the real axis, introduce the transformation 
\begin{equation*}
  X^{(2)}(z)=
  \begin{cases}
   X^{(1)}(z)\widehat{w}^{\frac{\sigma_3}{2}}(z)\begin{pmatrix}
    0& 1\\
    -1& 0
\end{pmatrix}^{-1}, & \text{for $z$ in regions $\widetilde{\Omega}_{IV}$, $\widetilde{\Omega}_{V.1}$ , $\widetilde{\Omega}_{V.2}$, $\widetilde{\Omega}_{VI}$},\\
X^{(1)}(z)\widehat{w}^{\frac{\sigma_3}{2}}(z), & \text{otherwise.}
\end{cases}
\end{equation*}
This yields the jump relations for $X^{(2)}(z)$ given in Figure \ref{Contour X2}.

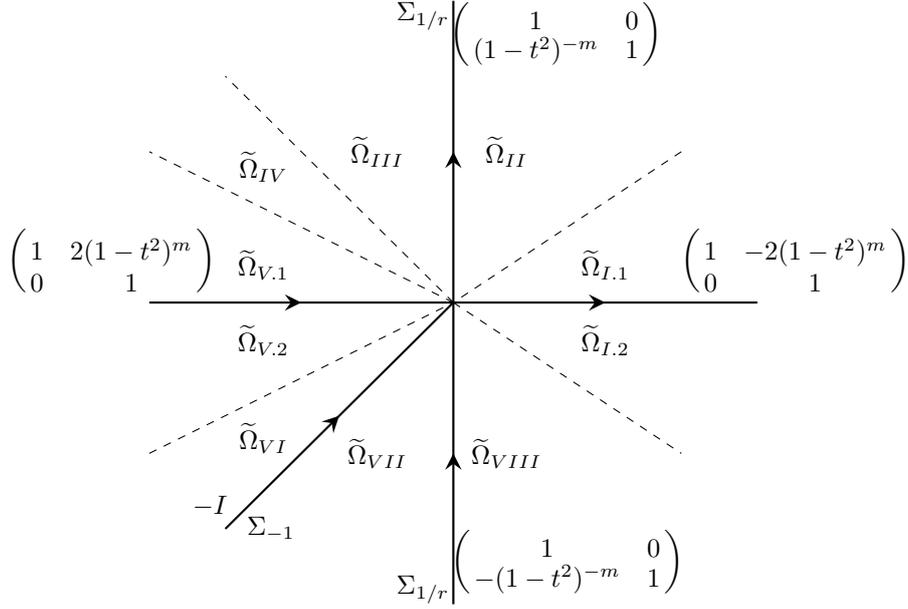
\begin{figure}[H]
    \centering
    \begin{tikzpicture}[
        myarrow/.style={
            thick,
            postaction={decorate,
            decoration={markings,mark=at position #1 with {\arrow[scale=1.5]{stealth}}}}
        },
        every matrix/.style={scale=0.4} 
    ]
        \draw[dashed] (0,0) -- (3,-2);
        \draw[dashed] (3,2) -- (0,0);

         \draw[myarrow=0.5] (0,0) -- (4,0);
        \draw[myarrow=0.5] (-4,0) -- (0,0);
        
        \draw[myarrow=0.5] (0,0) -- (0,4);
        \draw[dashed] (0,0) -- (-3,3);
        \draw[dashed] (0,0) -- (-4,2);

        \draw[dashed] (-4,-2) -- (0,0) ;
        \draw[myarrow=0.5] (-3,-3) -- (0,0) ;
        \draw[myarrow=0.5] (0,-4) -- (0,0) ;
        
        \matrix (m4) at (1.35,3.5) [matrix of math nodes, left delimiter=(, right delimiter=), inner sep=0pt, column sep=10pt, row sep=3pt] {
            1 & 0\\
            (1-t^2)^{-m}& 1 \\
        };

        \matrix (m4) at (1.5,-3.5) [matrix of math nodes, left delimiter=(, right delimiter=), inner sep=0pt, column sep=10pt, row sep=3pt] {
            1 & 0\\
            -(1-t^2)^{-m} & 1 \\
        };
        \matrix (m4) at (4.5,0.5) [matrix of math nodes, left delimiter=(, right delimiter=), inner sep=0pt, column sep=10pt, row sep=3pt] {
            1 &-2(1-t^2)^m\\
            0 & 1 \\
        };
        \matrix (m4) at (-4.5,0.5) [matrix of math nodes, left delimiter=(, right delimiter=), inner sep=0pt, column sep=10pt, row sep=3pt] {
            1 &2(1-t^2)^m\\
            0 & 1 \\
        };
        \node (eq) at (-3.2,-2.7) {$-I$}; 
           
        \node (eq) at (2,0.5) {$\widetilde{\Omega}_{I.1}$};
        \node (eq) at (2,-0.5) {$\widetilde{\Omega}_{I.2}$};
        \node (eq) at (0.7,2) {$\widetilde{\Omega}_{II}$};
        \node (eq) at (0.7,-2) {$\widetilde{\Omega}_{VIII}$};
        \node (eq) at (-1,2) {$\widetilde{\Omega}_{III}$};
        \node (eq) at (-1,-2) {$\widetilde{\Omega}_{VII}$};
        \node (eq) at (-2.5,1.8) {$\widetilde{\Omega}_{IV}$};
        \node (eq) at (-2.5,-1.8) {$\widetilde{\Omega}_{VI}$};
        \node (eq) at (-2.5,0.5) {$\widetilde{\Omega}_{V.1}$};
        \node (eq) at (-2.5,-0.5) {$\widetilde{\Omega}_{V.2}$}; 

        \node (eq) at (-2.4,-3) {$\Sigma_{-1}$};
        \node (eq) at (-0.4,-3.8) {$\Sigma_{1/r}$};
        \node (eq) at (-0.4,3.8) {$\Sigma_{1/r}$};
    \end{tikzpicture}
    \caption{Jump matrices for $X^{(3)}(z)$.}
    \label{Contour X3}
\end{figure}

\noindent
\textbf{Second transformation:} To collapse the jumps on $\Sigma_r\cap u_{\varepsilon}$ and $\Sigma_1\cap u_{\varepsilon}$ to the real axis, define 
\begin{equation*}
X^{(3)}(z)=\begin{cases}
    X^{(2)}(z)\begin{pmatrix}
        1& -(1-z^2)^{m}\\
        0 & 1
    \end{pmatrix}^{-1},  & \text{for $z$ in regions $\widetilde{\Omega}_{I.2}$ and $\widetilde{\Omega}_{V.1}$},\\
X^{(2)}(z)\begin{pmatrix}
        1& (1-z^2)^{m}\\
        0 & 1
    \end{pmatrix}^{-1}, & \text{for $z$ in regions $\widetilde{\Omega}_{I.1}$ and $\widetilde{\Omega}_{V.2}$},\\
    X^{(2)}(z),& \text{otherwise}.  
\end{cases}   
\end{equation*}
This leads to the jump relations for $X^{(3)}(z)$ shown in Figure \ref{Contour X3}.

The RHP for $X^{(3)}(z)$ is solved using the Parabolic Cylinder RHP \eqref{RHP parabolic}, by introducing a change of variables near origin. Indeed, note that for $|z|$ small, we can write 
$$(1-z^2)^m=e^{\phi(z)},$$
where $\phi(z)=m\log(1-z^2)$. The function $\phi$ has a critical point at $z=0$, with $\phi''(0)=-2m\neq 0$. Hence, there exist an open neighbourhood $u_\varepsilon$ of $z=0$ and a biholomorphic map 
\begin{equation}\label{Eq def map lambda}
\begin{split}
    \lambda:\,\, &u_{\varepsilon}\rightarrow U_m(0)\\
             &z\rightarrow \xi:=\lambda(z),
\end{split}  
\end{equation}
such that
\begin{equation*}
   \phi(\xi)=\phi(\lambda(z))=-\frac{i\xi^2}{2}, \hspace{0.5cm} \text{$z\in u_{\varepsilon}$}.
\end{equation*}
More explicitly, the map $\lambda$ can be written as
\begin{equation}\label{Eq def xi(z)}
  \xi=\lambda(z)=\sqrt{2m}\,z\,e^{-\pi i/4}\,\tau(z),
\end{equation}
 where
 $$\tau(z)=\sqrt{\sum_{k=0}^{\infty}\frac{z^{2k}}{k+1}}.$$ This last function is analytic in the disk $u_\varepsilon$. For $z\in u_{\varepsilon}$ 
 \begin{equation}\label{Eq identity z and xi}
     (1-z^2)^m=e^{-\frac{i\xi^2}{2}}.
 \end{equation}
Under this change of variables the segments $(-\varepsilon, \varepsilon)$ and $(-i\varepsilon, i\varepsilon)$ are rotated by an angle of $\pi/4$ in the clockwise direction. The part of the curve $\Sigma_{-1}$ in the third quadrant of Figure \ref{Contour X3} is mapped to a straight line segment in the negative real axis of the $\xi$ plane. This gives the jump relations in the $\xi$-plane shown in Figure \ref{Contour X3 xi}:
\begin{figure}[H]
    \centering
    \begin{tikzpicture}[
        myarrow/.style={
            thick,
            postaction={decorate,
            decoration={markings,mark=at position #1 with {\arrow[scale=1.5]{stealth}}}}
        },
        every matrix/.style={scale=0.4} 
    ]
        \draw[myarrow=0.7] (-3,-3) -- (3,3);
        \draw[myarrow=0.3] (-3,-3) -- (3,3);
        \draw[myarrow=0.3] (-5,0) -- (0,0);
        \draw[myarrow=0.3] (-3,3) -- (3,-3);
        \draw[myarrow=0.7] (-3,3) -- (3,-3);
        
        \matrix (m1) at (3.5,2) [matrix of math nodes, left delimiter=(, right delimiter=), inner sep=0pt, column sep=10pt, row sep=3pt] {
    1 & 0 \\
    a^2 & 1 \\
};
        
        \matrix (m2) at (3.5,-2) [matrix of math nodes, left delimiter=(, right delimiter=), inner sep=0pt, column sep=10pt,row sep=3pt] {
            1 & -2e^{-\frac{i\xi^2}{2}}  \\
           0& 1\\
        };
        \matrix (m3) at (-4,2) [matrix of math nodes, left delimiter=(, right delimiter=), inner sep=0pt, column sep=10pt, row sep=3pt] {
            1 & 2 e^{-\frac{i\xi^2}{2}}\\
           0 & 1 \\
        };
        \matrix (m4) at (-4,-2) [matrix of math nodes, left delimiter=(, right delimiter=), inner sep=0pt, column sep=10pt, row sep=3pt] {
            1 & 0\\
            -e^{\frac{i\xi^2}{2}} & 1 \\
        };
        \node (eq) at (-5,0.3) {$-I$};
    \end{tikzpicture}
    \caption{Jump matrices for $X^{(3)}(\xi)$.}
    \label{Contour X3 xi}
\end{figure}
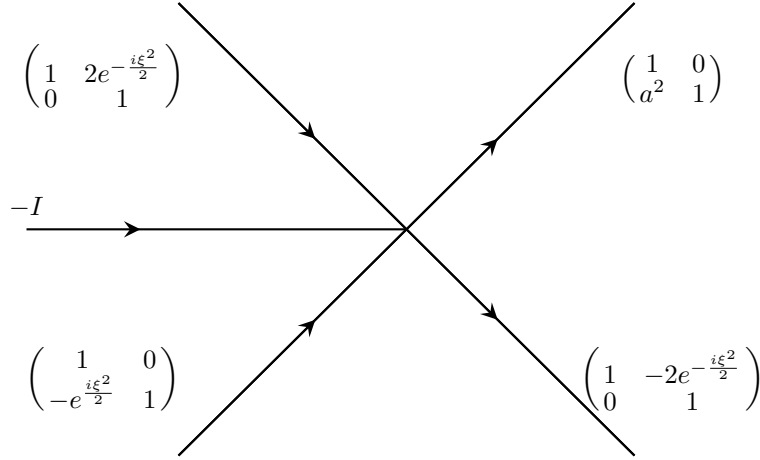

\noindent
\textbf{Third transformation:} In order to eliminate the exponential terms in the jumps for $X^{(3)}$, and to adjust the coefficients in the off diagonal entries of the jump matrices, introduce the transformation:
\begin{equation*}
    X^{(4)}(\xi)=a^{-\sigma_3}X^{(3)}(\xi)\,a^{\sigma_3} e^{-\frac{i\xi^2}{4}\sigma_3}.
\end{equation*}
The constant $a$ will be chosen in the next transformation. The jump relations for $X^{(4)}(\xi)$ are now shown in Figure \ref{Contour X4 xi}:
\begin{figure}[H]
    \centering
    \begin{tikzpicture}[
        myarrow/.style={
            thick,
            postaction={decorate,
            decoration={markings,mark=at position #1 with {\arrow[scale=1.5]{stealth}}}}
        },
        every matrix/.style={scale=0.5} 
    ]
        \draw[myarrow=0.7] (-3,-3) -- (3,3);
        \draw[myarrow=0.3] (-3,-3) -- (3,3);
        \draw[myarrow=0.3] (-5,0) -- (0,0);
        \draw[myarrow=0.3] (-3,3) -- (3,-3);
        \draw[myarrow=0.7] (-3,3) -- (3,-3);
        
        \matrix (m1) at (3.5,2) [matrix of math nodes, left delimiter=(, right delimiter=), inner sep=0pt, column sep=10pt, row sep=3pt] {
            1 & 0 \\
           a^2  & 1\\
        };
        
        \matrix (m2) at (3.7,-2) [matrix of math nodes, left delimiter=(, right delimiter=), inner sep=0pt, column sep=10pt, row sep=3pt] {
            1 & -2a^{-2} \\
           0& 1\\
        };
        \matrix (m3) at (-4,2) [matrix of math nodes, left delimiter=(, right delimiter=), inner sep=0pt, column sep=10pt, row sep=3pt] {
            1 & 2 a^{-2}\\
           0 & 1 \\
        };
        \matrix (m4) at (-4,-2) [matrix of math nodes, left delimiter=(, right delimiter=), inner sep=0pt, column sep=10pt, row sep=3pt] {
            1 & 0\\
            -a^2 & 1 \\
        };
        \node (eq) at (-5,0.3) {$-I$};
    \end{tikzpicture}
    \caption{Jump matrices for $X^{(4)}(\xi)$.}
    \label{Contour X4 xi}
\end{figure}
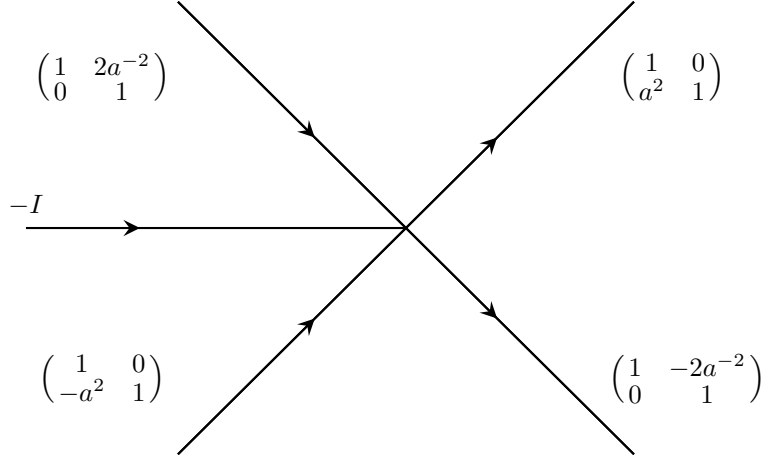

\noindent
\textbf{Fourth transformation:} The RHP for $X^{(4)}(\xi)$ can be transformed into RHP \ref{RHP parabolic} by  choosing $a=2^{1/4}e^{-\frac{i\pi}{4}}$ and introducing a rotation of $\frac{5\pi}{4}$ in the clockwise direction, that is, defining $\eta=\xi e^{-\frac{5\pi i}{4}}$. Indeed,
\begin{equation*}
     X^{(5)}(\eta)= X^{(4)}(\eta\,e^{\frac{5\pi i}{4}} ),
\end{equation*}
satisfies the jump relations shown in Figure  \ref{Contour X^{(5)}}.
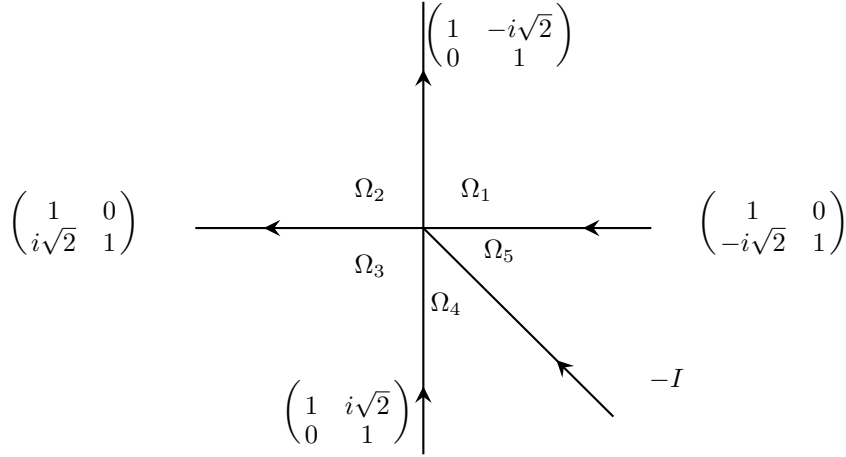
\begin{figure}[htbp]
    \centering
    \begin{tikzpicture}[
        myarrow/.style={
            thick,
            postaction={decorate,
            decoration={markings,mark=at position #1 with {\arrow[scale=1.5]{stealth}}}}
        },
        every matrix/.style={scale=0.5} 
    ]
        \draw[myarrow=0.7] (0,0) -- (0,3);
        \draw[myarrow=0.3] (0,-3)--(0,0);
        \draw[myarrow=0.7] (0,0) -- (-3,0);
        \draw[myarrow=0.3] (3,0)--(0,0);
        \draw[myarrow=0.3] (2.5,-2.5) -- (0,0);

        \node (eq) at (0.7,0.5) {$\Omega_1$};
        \node (eq) at (-0.7,0.5) {$\Omega_2$};
        \node (eq) at (0.3,-1) {$\Omega_4$};
        \node (eq) at (-0.7,-0.5) {$\Omega_3$};
        \node (eq) at (1,-0.3) {$\Omega_5$};

        \node (eq) at (3.2,-2) {$-I$};
        
        \matrix (m2) at (4.6,0) [matrix of math nodes, left delimiter=(, right delimiter=), inner sep=0pt, column sep=10pt, row sep=3pt] {
            1 & 0  \\
           -i\sqrt{2}& 1\\
        };
        
        \matrix (m3) at (1,2.5) [matrix of math nodes, left delimiter=(, right delimiter=), inner sep=0pt, column sep=10pt, row sep=3pt] {
            1 &  -i\sqrt{2}\\
           0 & 1 \\
        };
        
        \matrix (m4) at (-4.6,0) [matrix of math nodes, left delimiter=(, right delimiter=), inner sep=0pt, column sep=10pt, row sep=3pt] {
            1 & 0\\
            i\sqrt{2}& 1 \\
        };

        \matrix (m3) at (-1,-2.5) [matrix of math nodes, left delimiter=(, right delimiter=), inner sep=0pt, column sep=10pt, row sep=3pt] {
            1 &  i\sqrt{2}\\
           0 & 1 \\
        };
        
    \end{tikzpicture}
    \caption{Contour $\Sigma_{X^{(5)}}$ and  jump matrices $J_{X^{(5)}}(\eta)$.}
    \label{Contour X^{(5)}}
\end{figure}
Observe that changing the orientation of the three rays pointing to the origin, this problem is equivalent to RHP \ref{RHP parabolic}, with $\nu=-\frac{1}{2}$. Hence, a solution for this problem is
\begin{equation*}
    X^{(5)}(\eta)=Z(\eta),
\end{equation*}
where $Z(\eta)$ is the solution of RHP \ref{RHP parabolic} defined in \eqref{Eq sol RHP for P}.

In the variable $\xi$, this yields
\begin{equation*}
    X^{(4)}(\xi)=X^{(5)}(\xi \, e^{-\frac{5\pi i}{4}})=Z(\xi \, e^{-\frac{5\pi i}{4}}).
\end{equation*}
Now, as $\xi\rightarrow \infty$, \eqref{Eq Asymptotics PC} gives
\begin{equation}\label{Eq Asympt Z}
          Z(\xi \, e^{-\frac{5\pi i}{4}})=\frac{(\xi \, e^{-\frac{5\pi i}{4}})^{-\frac{\sigma_3}{2}}}{\sqrt{2}} \begin{pmatrix}
            1+\frac{3i}{8\xi^2}+\mathcal{O}\left(\frac{1}{\xi^4}\right)&1-\frac{3i}{8\xi^2}+\mathcal{O}\left(\frac{1}{\xi^4}\right)\\
            1-\frac{5i}{8\xi^2}+\mathcal{O}\left(\frac{1}{\xi^4}\right)&-1-\frac{5i}{8\xi^2}+\mathcal{O}\left(\frac{1}{\xi^4}\right)
        \end{pmatrix} e^{-\frac{i\xi^2}{4}\sigma_3}.
\end{equation}
The matrix function $$\widetilde{U}_l(z):=X^{(1)}(\xi(z)),$$
satisfying the jump relations from the local RHP \ref{RHP local} is now obtained by unfolding all the transformations performed to produce $X^{(4)}(\xi)$.  However, an adjustment is required to satisfy the third condition, which concerns the jump relation on the boundary of $u_{\varepsilon}$ where the local transformation $\lambda$ is defined (see \eqref{Eq def map lambda}), this being of the form $I+\mathcal{O}(1/m)$. This is addressed in the next subsection.
\subsection{Global solution} 
\label{sec Global solution}

Subsection \ref{sec Outside Solution} gave the asymptotic solution $N(z)$, to the RHP for $U(z)$ away from the origin. Subsection \ref{sec local parametrix} gave a matrix function $\widetilde{U}_l(z)$, satisfying the same jump relations as $U(z)$ within a neighborhood of the origin. The next goal is to use these two solutions to obtain the uniform global asymptotics in the complex plane. To achieve this, the task is divided into two steps. In the first step, the local solution \(\widetilde{U}_l(z)\) is multiplied on the left by an auxiliary matrix \(Q(z)\) to obtain a jump of the form \(I + \text{small}\) on \(\Sigma_{\varepsilon} := \partial u_{\varepsilon}\). This ensures the local solution satisfies the third condition of RHP~\ref{RHP local}. However, this will introduce a pole at the origin, an issue addressed in the second step.

\subsubsection{\textbf{Adjustment of the local solution}}
\label{subsec Adjust local sol}
Summary up to now:
\begin{itemize}
    \item $Y(z)$: solution to the original problem on the contour $\Sigma.$
    \item $U(z)$: solution to the transformed problem with lenses opened around $\Sigma_{-1}$.
    \item $N(z)$: approximation for $U(z)$ away from the origin, ignoring negligible jumps ($I+\mathcal{O}(1/m)$), and  considering only the jump on $\Sigma_{-1}$.
    \item $\widetilde{U}_l(z)$: local approximation for $U(z)$ on $u_{\varepsilon}$, obtained by unravelling all the transformations leading to $X^{(4)}(\xi(z))=Z(\xi(z)e^{\frac{-5\pi i}{4}})$.    
\end{itemize}
Now define $$ U_l(z) = Q(z) \widetilde{U}_l(z), $$ where the pre-factor matrix $Q(z)$ is a meromorphic matrix-valued function to be chosen in \eqref{EQ Def Q}, so that the local solution meets the appropriate jump condition on the boundary of the disk $u_{\varepsilon}$. Consider the matrix defined by
\begin{equation}\label{EQ Error matrix for U}
    E(z)=\begin{cases}
        U(z)N^{-1}(z),       & z\notin u_{\varepsilon};\\
        U(z)\,U_l^{-1}(z), & z \in u_{\varepsilon}.
    \end{cases}
\end{equation}
The jump relations for $E(z)$ are computed next. Outside the disc $u_{\varepsilon}$, one finds
\begin{equation*}
    J_E=(UN^{-1})_-^{-1}(UN^{-1})_+=N_-J_UN_+^{-1}.
\end{equation*}
Recall that $ N(z)$  satisfies the same jump relation as $U(z) $ on $\Sigma_{-1} $ and is analytic elsewhere (outside of $u_{\varepsilon}$).
Using the jump relations for $U(z)$ (see \eqref{Eq jump U}), and the definition of $N(z)$ (see \eqref{Eq N explicit}), it follows that
\begin{equation}\label{EQ jump for E out}
    J_E(t)=\begin{cases}
        I, & \text{ $t\in \Sigma_{-1}$},\\
        \begin{pmatrix}
       1& -\frac{(1-t^2)^{m}\left(\frac{1-t}{t}\right)}{\widetilde{w}(t)}\\
      0 & 1
    \end{pmatrix}, & \text{$t\in \Sigma_{r}$},\\
        \begin{pmatrix}
       1& 0,\\
      \frac{\left(\frac{t}{1+t}\right)}{(1-t^2)^{m}\widetilde{w}(t)} & 1
    \end{pmatrix}, & \text{$t\in \Sigma_{1/r}$},\\
        \begin{pmatrix}
       1 & (1-t^2)^m\widetilde{w}(t)\left(\frac{1+t}{t}\right),\\
       0 & 1 \\
    \end{pmatrix}, & \text{ $t\in \Sigma_{1}$.}
    \end{cases}
\end{equation}
Therefore $E(z)$ has no jump on $\Sigma_{-1}$ outside $u_{\varepsilon}$, and it has negligible jumps ($I+\mathcal{O}(1/m)$) on  $\Sigma_{r}$, $\Sigma_{1/r}$, and $\Sigma_{1}$. 

On the other hand, the error matrix $E(z)$ is analytic in the entire neighborhood $u_{\varepsilon}$, since $U_l$ satisfies the same jump relations as $U(z)$. Finally, $E(z)$ has also jump relations on $\Sigma_{\varepsilon}$. These should be of the form $I + small$ (see RHP \eqref{RHP local}). These jumps are given by
\begin{equation}\label{EQ jump E uepsilon}
    J_{E_{\varepsilon}}=E_-^{-1}E_+=[(UN^{-1})^{-1}]_-[UU_l^{-1}]_+=NU_-^{-1}U_+U_l^{-1}=NU_l^{-1},
\end{equation}
since $U(z)$ has no jumps on the boundary of $u_{\varepsilon}$. The jump relations for $E(z)$ on $\Sigma_{\varepsilon}$ must be analyzed independently in each of the regions of Figure \ref{Contour X1}. As $m\rightarrow\infty$ and for $z\in \Sigma_{\epsilon}$, the values $U_l(z)$  are related to the asymptotics of $Z(\xi \, e^{-\frac{5\pi i}{4}})$, as $\xi\rightarrow \infty$. In detail

\begin{equation}\label{Eq Asympt ZZ}
          Z(\xi \, e^{-\frac{5\pi i}{4}})\,e^{\frac{i\xi^2}{4}\sigma_3} =\frac{1}{\sqrt{2}}(\xi \, e^{-\frac{5\pi i}{4}})^{-\frac{\sigma_3}{2}}Z^{(0)}\left[
          I+\frac{1}{\xi^2}\left(Z^{(0)}\right)^{-1}Z^{(1)}+\mathcal{O}\left(\frac{1}{\xi^4}\right)\right].
\end{equation}
Here
\begin{equation*}
    Z^{(0)}=\begin{pmatrix}
            1 &1\\
            1&-1
        \end{pmatrix} \hspace{0.5cm}\text{and}\hspace{0.5cm} Z^{(1)}=\begin{pmatrix}
                 \frac{3i}{8} & -\frac{3i}{8} \\[2pt]
                -\frac{5i}{8} & -\frac{5i}{8}
                 \end{pmatrix}.
\end{equation*}
Consider $z$ in region $\widetilde{\Omega}_{III}$, then
\begin{align*} 
U_l(z)=&\,\,Q(z)\,\widetilde{U}_l(z)\\
=&\,\,Q(z)\,X^{(1)}(z)\\
=&\,\,Q(z)\,X^{(3)}(z)\,\widehat{w}(z)^{-\frac{\sigma_3}{2}}\\
=&\,\,Q(z)\,a^{\sigma_3}X^{(4)}(z)\,e^{\frac{i\xi^2}{4}\sigma_3}\,a^{-\sigma_3}\widehat{w}(z)^{-\frac{\sigma_3}{2}}\\
=&\,\,Q(z)\,a^{\sigma_3}Z\left(\xi(z) \, e^{\frac{-5\pi i}{4}}\right)e^{\frac{i\xi^2}{4}\sigma_3}\,a^{-\sigma_3}\widehat{w}(z)^{-\frac{\sigma_3}{2}}.
\end{align*}
Using the asymptotics for $Z(\xi \, e^{-\frac{5\pi i}{4}})$ given in  \eqref{Eq Asympt ZZ}, the jump relation \eqref{EQ jump E uepsilon} (in region $\widetilde{\Omega}_{III}$) is given by
\begin{align*}
   J_{E_{\varepsilon}}(z)&=N(z)U_l^{-1}(z)\\
       &=N(z) \widehat{w}^{\frac{\sigma_3}{2}}(z) a^{\sigma_3} (Z(\xi \, e^{-\frac{5\pi i}{4}})e^{\frac{i\xi^2}{4}\sigma_3})^{-1} a^{-\sigma_3} Q^{-1}(z)\\
       &=\sqrt{2}N(z) \widehat{w}(z)^{\frac{\sigma_3}{2}} a^{\sigma_3} \left[
          I-\frac{1}{\xi^2}\left(Z^{(0)}\right)^{-1}Z^{(1)}+\mathcal{O}\left(\frac{1}{\xi^4}\right)\right]\\
          &\hspace{4.5cm}\cdot\left(Z^{(0)}\right)^{-1}\left(\xi \, e^{-\frac{5\pi i}{4}}\right)^{\frac{\sigma_3}{2}}a^{-\sigma_3} Q^{-1}(z).
\end{align*}
This suggests the introduction of 
\begin{equation}\label{EQ Def Q}
    Q(z)=\sqrt{2}\, N(z) \widehat{w}(z)^{\frac{\sigma_3}{2}} a^{\sigma_3}\left(Z^{(0)}\right)^{-1}\left(\xi \, e^{-\frac{5\pi i}{4}}\right)^{\frac{\sigma_3}{2}}a^{-\sigma_3},
\end{equation}
since it leads to the jump on $\Sigma_{\varepsilon}$ having the form $I+\mathcal{O}(1/m)$, which is negligible  as $m\rightarrow\infty$. It can be shown that the same $Q(z)$ yields negligible jumps for $E(z)$ over all $\Sigma_{\varepsilon}$.

\begin{figure}[htbp]
    \centering  
\begin{tikzpicture}[scale=0.8]
    \node at (-1,0) [circle,fill,inner sep=1pt,label=below:-1] {};
    \node at (1,0) [circle,fill,inner sep=1pt,label=below:1] {};

    \draw[line width=1pt]
    plot[domain=90:270,samples=400,smooth] ({sqrt(3.5)*cos(\x)},{sqrt(3)*sin(2*\x)});

     \draw[line width=1pt]
    plot[domain=-180:180,samples=200,smooth] ({sqrt(2)*cos(\x)},{0.6*sqrt(2)*sin(2*\x)});
    
    \draw[line width=1pt, blue!100!] (0,0) circle [radius=0.37];
    
    \fill[white] (0,0) circle [radius=0.35];

    \node (eq) at (2.8,0) {\small{$I+O(e^{-m})$}};
   \node (eq2) at (0.5,1.3) [blue!100!] {\small{$I+O(\frac{1}{m})$}};
    
\end{tikzpicture}
    \caption{Jump diagram for $E$ over the contour $\widetilde{\Sigma}$.}
    \label{Contour_sigma_hat}
\end{figure}
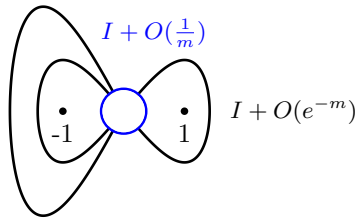

From here it follows that
\begin{equation}
\label{Eq Error E limit}
 J_{E}(z)=I+\mathcal{O}\left(\frac{1}{m}\right).
\end{equation}
This shows that the RHP for $U(z)$ has become an RHP for $E(z)$, with jumps that are uniformly $I+\mathcal{O}(1/m)$ on all contours (see Figure \ref{Contour_sigma_hat}).

However, the matrix $Q(z)$ introduces a pole at the origin, coming from the terms $N(z)$ and $\xi(z)^{\frac{\sigma_3}{2}}$. Indeed, from \eqref{EQ Def Q},
\begin{align*}
    Q(z)&=\sqrt{2} \left[\frac{(1+z)^{\frac{1}{2}}}{z^{\frac{1}{2}}}\right]^{\sigma_3} \widehat{w}(z)^{\frac{\sigma_3}{2}}a^{\sigma_3}\begin{pmatrix}
            1 &1\\
            1&-1
        \end{pmatrix}^{-1}\left(\xi \, e^{-\frac{5\pi i}{4}}\right)^{\frac{\sigma_3}{2}}a^{-\sigma_3}\\
        &=\frac{\sqrt{2}}{2}\left[(1+z)\widehat{w}(z)\right]^{\frac{\sigma_3}{2}}\left[z^{-\frac{1}{2}}a\right]^{\sigma_3}\begin{pmatrix}
            1 &1\\
            1&-1
        \end{pmatrix} \left[z^{-\frac{1}{2}}a\right]^{-\sigma_3}\left[\frac{\xi\,e^{-\frac{5\pi i}{4}}}{z}\right]^{\frac{\sigma_3}{2}}\\
        &=\frac{\sqrt{2}}{2}\left[(1+z)\widehat{w}(z)\right]^{\frac{\sigma_3}{2}}\begin{pmatrix}
           1 &\frac{a^2}{z}\\
            \frac{z}{a^2}&-1
        \end{pmatrix}\left[\frac{\xi\,e^{-\frac{5\pi i}{4}}}{z}\right]^{\frac{\sigma_3}{2}}\\
        &=\frac{\sqrt{2}}{2}Q_L(z) \begin{pmatrix}
           1 &\frac{a^2}{z}\\
            \frac{z}{a^2}&-1
        \end{pmatrix}   
        Q_R(z),
\end{align*}
where the matrix functions
\begin{equation}\label{Eq def Q_L}
    Q_L(z)=\left[(1+z)\widehat{w}(z)\right]^{\frac{\sigma_3}{2}},\hspace{0.5cm}\text{and}\hspace{0.5cm} Q_R(z)=\left[\frac{\xi\,e^{-\frac{5\pi i}{4}}}{z}\right]^{\frac{\sigma_3}{2}},
\end{equation}
are bounded near the origin in view of \eqref{Eq def xi(z)}. 

Let 
\begin{equation}\label{Eq. A matrices}
    A_0=\sigma_3, \hspace{0.5cm} 
    A_1=\begin{pmatrix}
        0&a^2\\
        0&0
    \end{pmatrix},\hspace{0.5cm}\text{and}\hspace{0.5cm} A_2=\begin{pmatrix}
        0&0\\
        a^{-2}&0
    \end{pmatrix}.
\end{equation}
Then $Q(z)$ can be written as:
\begin{equation}\label{Eq. Matrix Q}
Q(z)=\frac{\sqrt{2}}{2}Q_L(z)\left[A_0+\frac{1}{z}A_1+z\,A_2\right]Q_R(z),   
\end{equation}
which shows the appearance of a pole in the middle term. 

This implies that the matrix $Q$ introduces a pole in the local parametrix $U_l$. However, $U_l$ must be analytic at zero. The issue introduced by this pole is addressed in the next subsection.

\subsubsection{\textbf{Addressing the issue of the pole introduced at $z=0$.}} 
Let $\widetilde{E}(z)$ be a matrix function with the same jump conditions as $E(z)$ (see Figure \ref{Contour_sigma_hat}), but now require $\widetilde{E}(z)$ to be analytic at zero. This function can be constructed using Neumann series, since all its jump matrices are close to the identity matrix for $m$ large. Below it is shown that $\widetilde{E}(z)=I+\mathcal{O}(1/m)$ uniformly for $z$ bounded away from the contour $\widetilde{\Sigma}$. Indeed, the RHP for $\widetilde{E}$ is:
\begin{equation*}
\begin{cases}
\widetilde{E}(z)\,\,\,\text{is analytic in}\,\, \mathbb{C}\setminus \widetilde{\Sigma},\\
\widetilde{E}_+(t)=\widetilde{E}_-(t)J_{E}(t),\,\,\,\text{for}\,\,z\in \widetilde{\Sigma}.\\
\widetilde{E}(z)=I+\mathcal{O}\left(\frac{1}{z}\right), \,\,\,\text{when}\,\,z\rightarrow\infty.
\end{cases}    
\end{equation*}
The relation $\widetilde{E}_+=\widetilde{E}_-J_{E}=\widetilde{E}_-(I+J_{E}-I)$ can be rewritten as 
\begin{equation}\label{RHP Etilde-I}
    (\widetilde{E}-I)_+-(\widetilde{E}-I)_-=\widetilde{E}_-(J_{E}-I).
\end{equation}
Consider the Cauchy operator $\mathcal{C}$ on  $C(\widetilde{\Sigma})$, that is 
\begin{equation*}
    \mathcal{C}(F)(z)=\frac{1}{2\pi i}\oint\limits_{\widetilde{\Sigma}}\frac{F(t)}{t-z}\,dt,
\end{equation*}
with the corresponding Cauchy projection operators on the $``-"$ and $``+"$ sides of the contour:
\begin{equation*}
    \mathcal{C}_-(F)(\lambda)=\lim_{\begin{smallmatrix} z \to \lambda & \\ z\in - \text{side} \end{smallmatrix}} \mathcal{C}(F)(z),\hspace{1.5cm}\mathcal{C}_+(F)(\lambda)=\lim_{\begin{smallmatrix} z \to \lambda & \\ z\in + \text{side} \end{smallmatrix}} \mathcal{C}(F)(z),
\end{equation*}
and define the operator $\mathcal{C}_J$ on $L^2(\widetilde{\Sigma})$ by 
\begin{equation*}
    \mathcal{C}_J(F)(\lambda)= \mathcal{C}_-[F(J_{E}-I)](\lambda).
\end{equation*}
Using the Sokhotski Plemelj formula (See \cite{ablowitz2003complex}[Lemma 7.2.1]), equation \eqref{RHP Etilde-I} is seen to be equivalent to
\begin{equation}\label{Eq Etilde Cauchy 1}
    \widetilde{E}-I=\mathcal{C}(\widetilde{E}_-(J_{{E}}-I)).
\end{equation}
Then, the $``-"$ boundary value in \eqref{Eq Etilde Cauchy 1} implies that $\widetilde{E}_-$ satisfies the singular integral equation
\begin{equation}\label{Eq Singular Eq}
    (I-\mathcal{C}_{J})(\widetilde{E}_-)=I.
\end{equation}
It is shown next that there exists a constant 
$0<\kappa<\infty$, such that for $m$ large enough $||\mathcal{C}_{J}||_{L^2}\le\frac{\kappa}{m}$, so \eqref{Eq Singular Eq} can be solved by a Neumann series, that is 
\begin{equation*}
    \widetilde{E}_-=\sum_{k=0}^{\infty} \mathcal{C}_{J}^k(I).
\end{equation*}

Since the operator $\mathcal{C}_-$ is bounded in $L^2(\widetilde{\Sigma})$, there exists a constant $0<\widetilde{\kappa}<\infty$, such that for a matrix norm $\norm{\,\cdot\,}$  we have $||\norm{\mathcal{C}_-(F)}||_{L^2}\le \widetilde{\kappa} ||\norm{F}||_{L^2}$, for any matrix-valued function $F \in L^2(\widetilde{\Sigma})$. Then
\begin{equation}\label{Eq Bound1 CJ}
\begin{split}
    ||\norm{\mathcal{C}_{J}(F)}||_{L^2}&=||
    \norm{\mathcal{C}_{-}[F(J_{E}-I)]}||_{L^2}\\
    &\le \widetilde{\kappa}||\norm{F(J_{E}-I)}||_{L^2}\\
     &\le \widetilde{\kappa}||\norm{F}\norm{J_{E}-I}||_{L^2}\\
     &\le \widetilde{\kappa}||\norm{F}||\norm{J_{E}-I}||_{L^{\infty}}||_{L^2}\\
    &= \widetilde{\kappa}||\norm{F}||_{L^2}||\norm{J_{E}-I}||_{L^{\infty}}.
\end{split}
\end{equation}
In addition, \eqref{Eq Error E limit} shows the existence of a constant $\hat{\kappa}>0$ and $m_0\in \mathbb{N}$ such that for all $m\ge m_0$, 
\begin{equation}\label{Eq Bound JE-I}
    ||\norm{J_{E}-I}||_{L^{\infty}}\le \frac{\hat{\kappa}}{m}.
\end{equation}
Inequalities \eqref{Eq Bound1 CJ} and \eqref{Eq Bound JE-I} imply that, for $m\ge m_0$, $||\mathcal{C}_{J}||_{L^2}\le\frac{\kappa}{m}$, with $\kappa=\widetilde{\kappa}\hat{\kappa}$.

The bound on the norm of $\mathcal{C}_{J}$ shows the existence of $\widetilde{E}_-$ for $m$ sufficiently large. Moreover, it implies that
\begin{equation}\label{Eq Bound Etilde-I L2}
    \norm{\parallel\widetilde{E}_{-}\parallel}_{L^2}= ||\parallel I\parallel||_{L^2}+\mathcal{O}\left(\frac{1}{m}\right).
\end{equation}
Recall from \eqref{Eq Etilde Cauchy 1} that 
\begin{equation*} \label{Eq Etilde Cauchy 2}
    \widetilde{E}-I=\frac{1}{2\pi i}\oint\limits_{\widetilde{\Sigma}}\frac{\widetilde{E}_-(J_{E}-I)}{t-z}\,dt,
\end{equation*}
thus, for $z$ bounded away from the contour $\widetilde{\Sigma}$, say for $z$ in the set $$\{z\in \mathbb{C}: |z-s|\ge \delta>0, \,\,\,\forall s\in \widetilde{\Sigma}\},$$ one obtains 
\begin{equation*} 
\begin{split}
    \left\|\parallel\widetilde{E}-I\parallel\right\|_{L^{\infty}} &\le \frac{1}{2\pi} \oint\limits_{\widetilde{\Sigma}} \left\|\frac{1}{|t-z|} \left\|\widetilde{E}_-(J_E-I)\right\|\right\|_{L^{\infty}} \, d|t|, \\
    &\le \frac{1}{2\pi \delta} \oint\limits_{\widetilde{\Sigma}} \left\|\widetilde{E}_-(J_E-I)\right\| \, d|t|, \\
    &\le \frac{1}{2\pi \delta} \oint\limits_{\widetilde{\Sigma}} \left\|\widetilde{E}_-\right\|||\left\|J_E-I\right\| ||_{L^{\infty}} \, d|t|, \\
    &\le \frac{\hat{\kappa}}{2\pi \delta m} \oint\limits_{\widetilde{\Sigma}} \left\|\widetilde{E}_-\right\| \, d|t|, \\
    &\le \frac{1}{2\pi \delta m } \left\|\parallel\widetilde{E}_-\parallel\right\|_{L^2}||\norm{I}||_{L^2},
\end{split}
\end{equation*}
using Cauchy–Schwarz in the inequality. 

Since $||\parallel\widetilde{E}_-\parallel ||_{L^2}$ is bounded (see \eqref{Eq Bound Etilde-I L2}), it follows that
\begin{equation}\label{Eq Etilde Cauchy bound}
\widetilde{E}=I+\mathcal{O}\left(\frac{1}{m}\right).
\end{equation}
Define $\mathcal{E}(z) = E(z) \widetilde{E}^{-1}(z).$ Then $\mathcal{E}(z)$ has no jumps in the entire plane, but it has a simple pole at the origin, introduced in $E(z)$ by the factor $Q^{-1}$ in the inverse of the local parametrix $U_l$, see \eqref{EQ Error matrix for U}. Then
\begin{equation}\label{Error matrix for E}
    \mathcal{E}(z) = E(z) \widetilde{E}^{-1}(z)=I+\frac{\mathcal{E}^{(1)} }{z}.
\end{equation}
The constant matrix $\mathcal{E}^{(1)}$ is determined as follows: \eqref{EQ Error matrix for U} and  \eqref{Error matrix for E} imply that
\begin{equation}\label{U full asymptotics}
    U(z)=\begin{cases}
    \left(I+\frac{\mathcal{E}^{(1)} }{z}\right)\widetilde{E}(z)\, N(z), & \text{for $z$ outside} \,\,u_{\varepsilon},\\
    \left(I+\frac{\mathcal{E}^{(1)} }{z}\right)\widetilde{E}(z) \,U_l(z), & \text{for $z$ inside}\,\, u_{\varepsilon}.
\end{cases}
\end{equation}
Since $U(z)$ does not have a pole at $z=0$, $\mathcal{E}^{(1)}$ can be determined explicitly using the restriction that $U(z)$ stays bounded at the origin. To achieve this, consider the representation of $U(z)$ inside $u_{\varepsilon}$ and in the region $\widetilde{\Omega}_{III}$:
\begin{align*}
    U(z)&=\left(I+\frac{\mathcal{E}^{(1)} }{z}\right)\widetilde{E}(z)\,U_l(z)\\
    &=\left(I+\frac{\mathcal{E}^{(1)} }{z}\right)\widetilde{E}(z)\, Q(z) \,X^{(1)}(z)\\
    &=\left(I+\frac{\mathcal{E}^{(1)} }{z}\right)\widetilde{E}(z)\frac{\sqrt{2}}{2}Q_L(z)\left[A_0+\frac{1}{z}A_1+z\,A_2\right]Q_R(z)X^{(1)}(z),
\end{align*}
with $$X^{(1)}(z)=\,a^{\sigma_3}Z\left(\xi(z) \, e^{\frac{-5\pi i}{4}}\right)e^{\frac{i\xi^2}{4}\sigma_3}a^{-\sigma_3}\widehat{w}(z)^{-\frac{\sigma_3}{2}}.$$
Then
\begin{equation}\label{EQ U(z) for finding residue}
\begin{split}
   \sqrt{2}\,z\,U(z) &= \left(zI + \mathcal{E}^{(1)}\right)\widetilde{E}(z)Q_L(z) 
    \left[A_0+\frac{1}{z}A_1+z\,A_2\right]
    Q_R(z)X^{(1)}(z)\\ 
    &=\frac{1}{z}C_{-1}+C_0+z\,C_1+z^2C_2.
    \end{split}
\end{equation}
with 
\begin{align*}    C_{-1}=&\mathcal{E}^{(1)}\,\widetilde{E}\,Q_L\,A_1\,Q_R\,X^{(1)},\\    C_0=&\widetilde{E}\,Q_L\,A_1\,Q_R\,X^{(1)}+\mathcal{E}^{(1)}\,\widetilde{E}\,Q_L\,A_0\,Q_R\,X^{(1)}, \\    C_1=&\widetilde{E}\,Q_L\,A_0\,Q_R\,X^{(1)}+\mathcal{E}^{(1)}\widetilde{E}\,Q_L\,A_2\,Q_R\,X^{(1)}, \\
    C_2=& \widetilde{E}\,Q_L\,A_2\,Q_R\,X^{(1)}.
\end{align*}
Here, the matrix functions $C_j(z)$, $j=-1,0,1,2$, are analytic at the origin. Then 
$$C_{j}(z)=C_{j}^{(0)}+z\,C_{j}^{(1)}+z^2\,C_{j}^{(2)}+\cdots.$$
Since $U(z)$ is bounded at the origin,
$$\lim_{z\rightarrow 0} zU(z)=0.$$
Hence, \eqref{EQ U(z) for finding residue} yields
\begin{equation}\label{Eq. zU 1/z term}
    C_{-1}^{(0)}=C_{-1}(0)=0,
\end{equation}
and 
\begin{equation}\label{Eq. zU constant term}
 C_{-1}^{(1)}+C_{0}^{(0)}=C_{-1}'(0)+C_{0}(0)=0.
\end{equation}
The previous equations will be used to determine the constant matrix $\mathcal{E}^{(1)}$. First recall that $Q_R(0)$ and $X^{(1)}(0)$ are invertible. Then \eqref{Eq. zU 1/z term} yields
\begin{equation} \label{Eq. zU 1/z term simplified}
    \mathcal{E}^{(1)\,}\widetilde{E}(0)\,Q_L(0)\,A_1=0.
\end{equation}
Also, from the definition of $Q_L$ in \eqref{Eq def Q_L} it follows that
\begin{equation}\label{Eq. QL(0)=0}
    Q_L'(0)=0.
\end{equation}

Since
\begin{align*}
C_{-1}' = & \mathcal{E}^{(1)} \widetilde{E}' Q_L A_1 Q_R X^{(1)} + \mathcal{E}^{(1)} \widetilde{E} Q_L' A_1 Q_R X^{(1)}\\
&\hspace{1cm} + \mathcal{E}^{(1)} \widetilde{E} Q_L A_1 Q_R' X^{(1)}
+ \mathcal{E}^{(1)} \widetilde{E} Q_L A_1 Q_R X^{(1)'},
\end{align*}
equations \eqref{Eq. zU 1/z term simplified} and \eqref{Eq. QL(0)=0} yield
\begin{equation*}
   C_{-1}'(0) = \mathcal{E}^{(1)} \widetilde{E}'(0) Q_L(0) A_1 Q_R(0) X^{(1)}(0). 
\end{equation*}
Now, \eqref{Eq. zU constant term} gives
\begin{equation*}
 \mathcal{E}^{(1)} \widetilde{E}'(0)\, Q_L(0) A_1+\mathcal{E}^{(1)}\widetilde{E}(0)Q_L(0)A_0+ 
   \widetilde{E}(0)Q_L(0)A_1
   =0 , 
\end{equation*}
which implies that
\begin{equation*}
 \mathcal{E}^{(1)}=- \widetilde{E}(0)Q_L(0)A_1\,A_0^{-1}[Q_L(0)]^{-1}\,[\widetilde{E}(0)]^{-1}\left(I+ \widetilde{E}'(0)\, Q_L(0) A_1A_0^{-1}[Q_L(0)]^{-1}\,[\widetilde{E}(0)]^{-1}\right)^{-1}.   
\end{equation*}
The invertibility of the last term is guaranteed by  the fact that $\widetilde{E}'(0)=\mathcal{O}(1/m).$

Since $[Q_L(0)]^2=-1$, (see \eqref{Eq def Q_L}), and using the definitions of $A_0$ and $A_1$  (see \eqref{Eq. A matrices}), it follows that
\begin{equation*}
\begin{split}
    \mathcal{E}^{(1)}&=-\widetilde{E}(0)Q_L(0) \begin{pmatrix}
        0 & a^2 \\
        0 & 0
    \end{pmatrix}
    \begin{pmatrix}
        1&0\\
        0&-1
    \end{pmatrix}^{-1}
    [Q_L(0)]^{-1}[\widetilde{E}(0)]^{-1}\\
   &\hspace{4cm}\left(I+ \widetilde{E}'(0)\, Q_L(0) A_1A_0^{-1}[Q_L(0)]^{-1}\,[\widetilde{E}(0)]^{-1}\right)^{-1}\\
   &= -\widetilde{E}(0)\begin{pmatrix}
        0 & -a^2 \\
        0 & 0
    \end{pmatrix}[\widetilde{E}(0)]^{-1}\left(I+ \widetilde{E}'(0)\, Q_L(0) A_1A_0^{-1}[Q_L(0)]^{-1}\,[\widetilde{E}(0)]^{-1}\right)^{-1}.
\end{split}
\end{equation*}
Therefore, since $\widetilde{E}(0)=I+\mathcal{O} (1/m)$ ( see \eqref{Eq Etilde Cauchy bound}), one obtains
\begin{equation}\label{EQ Error 1}
    \mathcal{E}^{(1)}=\begin{pmatrix}
        0 & -a^2 \\
        0 & 0
    \end{pmatrix}+\mathcal{O}\left(\frac{1}{m}\right).
\end{equation}
It can be checked that this value of $\mathcal{E}^{(1)}$ is consistent with \eqref{Eq. zU 1/z term}.


Now all the necessary components to derive a global solution for $Y(z)$ across the entire complex plane have been established. Indeed, using \eqref{U full asymptotics}, and unfolding the initial transformations ($Y \rightarrow T \rightarrow U$), one obtains
\begin{equation}\label{Y full asymptotics}
    Y(z)=\begin{cases}
    \left(I+\frac{\mathcal{E}^{(1)} }{z}\right)\widetilde{E}(z)\, N(z)\,V(z), & \text{$z$ outside} \,\,u_{\varepsilon},\\
    \left(I+\frac{\mathcal{E}^{(1)} }{z}\right)\widetilde{E}(z) \,U_l(z)\,V(z), & \text{$z$ inside}\,\, u_{\varepsilon},
\end{cases}
\end{equation}
where
\begin{align*}
    V(z)&= K^{-1}(z)H^{-1}(z)\\
        &= \begin{pmatrix}
       1& 0\\
      -\frac{1}{(1-z^2)^{m}\,\widetilde{w}(z)}  & 1
    \end{pmatrix}
    \begin{pmatrix}
        (1+z)^m & 0 \\
        0 & (1+z)^{-m}
     \end{pmatrix} .
\end{align*}
Since the Jacobi polynomials $\pi_m(z)$ correspond to the entry $Y_{11}(z)$, it follows that
\begin{equation}\label{Eq Asymptotics for pi_m}
\begin{split}
     \pi_m(z)=&\left[\left(1+\frac{1}{z}
\mathcal{E}_{11}^{(1)}\right)\widetilde{E}_{11}+\frac{1}{z}\mathcal{E}_{12}^{(1)}\widetilde{E}_{21}\right]B_{11}\\
&\hspace{3cm}+\left[\left(1+\frac{1}{z}
\mathcal{E}_{11}^{(1)}\right)\widetilde{E}_{12}+\frac{1}{z}\mathcal{E}_{12}^{(1)}\widetilde{E}_{22}\right]B_{21},
\end{split}
\end{equation}
with
\begin{equation*}
    B(z)=\begin{cases}
    N(z)\,V(z), & \text{$z$ outside} \,\,u_{\varepsilon},\\
    U_l(z)\,V(z), & \text{$z$ inside}\,\, u_{\varepsilon}.
\end{cases}
\end{equation*}
In particular, for $z\in \Omega_-\setminus u_{\varepsilon}$,
\begin{align*}
    B(z)&= N(z)V(z)\\
        &= \begin{pmatrix}
        \frac{(1+z)^{\frac{1}{2}}}{z^{\frac{1}{2}}}& 0 \\
        0 & \frac{z^{\frac{1}{2}}}{(1+z)^{\frac{1}{2}}}
    \end{pmatrix}
    \begin{pmatrix}
       1& 0\\
      -\frac{1}{(1-z^2)^{m}\,\widetilde{w}(z)}  & 1
    \end{pmatrix}
    \begin{pmatrix}
        (1+z)^m & 0 \\
        0 & (1+z)^{-m}
     \end{pmatrix} \\
    &=  \begin{pmatrix}
       \frac{(1+z)^{\frac{1}{2}}}{z^{\frac{1}{2}}}  & 0 \\
        -\frac{z^{\frac{1}{2}}}{(1-z)^{\frac{1}{2}}(1-z^2)^m}&   \frac{z^{\frac{1}{2}}}{(1+z)^{\frac{1}{2}}}
    \end{pmatrix}
    \begin{pmatrix}
        (1+z)^m & 0 \\
        0 & (1+z)^{-m}
     \end{pmatrix}\\
     &= \begin{pmatrix}
       \frac{(1+z)^{\frac{1}{2}}(1+z)^m}{z^{\frac{1}{2}}}  & 0 \\
        -\frac{z^{\frac{1}{2}}}{(1-z)^{\frac{1}{2}}(1-z)^m}&   \frac{z^{\frac{1}{2}}}{(1+z)^{\frac{1}{2}}(1+z)^m}
    \end{pmatrix}.
\end{align*}
On the other hand,  for $z\in \Omega_-\cap u_{\varepsilon}$, using the definition of $U_l(z)$ in $\widetilde{\Omega}_{III}$, it follows that
\begin{equation}\label{B(z) near origin}
\begin{split}
B(z) &= U_l(z)V(z) \\
         &= Q(z)\widetilde{U}_l(z)V(z) \\
         &= \frac{\sqrt{2}}{2}\, N(z) \widehat{w}(z)^{\frac{\sigma_3}{2}} a^{\sigma_3}\begin{pmatrix}
            1 &1\\
            1&-1
        \end{pmatrix}\left(\xi \, e^{-\frac{5\pi i}{4}}\right)^{\frac{\sigma_3}{2}}Z\left(\xi\,e^{\frac{-5\pi i}{4}}\right)a^{-\sigma_3}\\
        &\hspace{7cm}\cdot e^{\frac{i\xi^2}{4}\sigma_3}\widehat{w}(z)^{-\frac{\sigma_3}{2}}V(z).
\end{split}
\end{equation}
Similar expressions for $B(z)$ can be derived in other regions of the complex plane. The region $\Omega_-$ is particularly important as it corresponds to the region where the zeros of the polynomials $P_m^{(m+1/2, -m-1/2)}(z)$ are located. This is discussed in the next section.
\section{Location of the zeros}
\label{sec zeros}
The monic Jacobi polynomials $\{\pi_m(z)\}$ appear in the upper-left corner of the matrix $Y(z)$ (see Proposition \ref{RHP Y sol polynomials}). In this section \eqref{Eq Asymptotics for pi_m} is used to determine the location of their zeros. The discussion is divided into two cases, depending on whether  $z$ is far from or close to the origin. This is required because the expressions for  $Y(z)$  differ in each case.
\subsection{Zeros away from the origin}

\begin{proof}[\textbf{Proof of Theorem \ref{Thm Asymptotics} (Part 1)}]
    From equation \eqref{Eq Asymptotics for pi_m}, the monic polynomials $\pi_m(z)$ are given in the region $ \Omega_-\setminus u_{\varepsilon}$ by
\begin{align*}    
\pi_m(z)&=B_{11}\left[\left(1+\frac{\mathcal{E}_{11}^{(1)}} 
    {z}\right)\widetilde{E}_{11}+ \frac{\mathcal{E}_{12}^{(1)}}{z}\widetilde{E}_{21}\right]+ B_{21}\left[\left(1+\frac{\mathcal{E}_{11}^{(1)}}{z}\right)\widetilde{E}_{12}+ \frac{\mathcal{E}_{12}^{(1)}}{z}\widetilde{E}_{22}\right]\\
    &=\frac{(1+z)^{\frac{1}{2}} (1+z)^m}{z^{\frac{1}{2}} }\left[\left(1+\frac{\mathcal{E}_{11}^{(1)}}{z}\right)\widetilde{E}_{11}+ \frac{\mathcal{E}_{12}^{(1)}}{z}\widetilde{E}_{21}\right]\\
    &  \hspace{3.5cm}-\frac{z^{\frac{1}{2}}}{(1-z)^{\frac{1}{2}}(1-z)^m}    \left[\left(1+\frac{\mathcal{E}_{11}^{(1)}}{z}\right)\widetilde{E}_{12}+ \frac{\mathcal{E}_{12}^{(1)}}{z}\widetilde{E}_{22}\right]\\
    &= \frac{(1+z)^{\frac{1}{2}} (1+z)^m}{z^{\frac{1}{2}}}\left\{\left(1+\frac{\mathcal{E}_{11}^{(1)}}{z}\right)\widetilde{E}_{11}+ \frac{\mathcal{E}_{12}^{(1)}}{z}\widetilde{E}_{21}\right.\\
    &\left. \hspace{1.5cm} -\frac{z}{(1-z)^{\frac{1}{2}}(1+z)^{\frac{1}{2}}(1-z^2)^m}\left[\left(1+\frac{\mathcal{E}_{11}^{(1)}}{z}\right)\widetilde{E}_{12}+ \frac{\mathcal{E}_{12}^{(1)}}{z}\widetilde{E}_{22}\right]\right\}.
\end{align*}

Recall from \eqref{Eq Etilde Cauchy bound} that $\widetilde{E}_{11}, \widetilde{E}_{22} = 1+\mathcal{O}(1/m)$ and $\widetilde{E}_{12}, \widetilde{E}_{21} = \mathcal{O}(1/m)$.  Also, \eqref{EQ Error 1} shows that $\mathcal{E}_{11}^{(1)}=\mathcal{O}(1/m)$ and  $\mathcal{E}_{12}^{(1)}=-a^2+\mathcal{O}(1/m)$. Thus, 
\begin{equation}\label{EQ pi_n(z)=F(q+r)}
    \pi_m(z)=\frac{(1+z)^{\frac{1}{2}} (1+z)^m}{z^{\frac{1}{2}} }\left\{1-\frac{i\sqrt{2}}{(1-z)^{\frac{1}{2}}(1+z)^{\frac{1}{2}}(1-z^2)^m}+r(m,z)\right\},
\end{equation}
 where $r(m,z)=\mathcal{O}(1/m)$, using $a^2=-i\sqrt{2}$.
\end{proof}
The next step is to show that, for sufficiently large $m$, the location of the zeros of $\pi_m(z)$ is determined from the solutions of
\begin{equation}\label{EQ Zeros prev}
    q(m,z):=1-\frac{i\sqrt{2}}{(1-z)^{\frac{1}{2}}(1+z)^{\frac{1}{2}}(1-z^2)^m}=0.    
\end{equation}
This is the statement of Theorem \ref{Thm Zeros out}.

\begin{figure}[H]
    \centering
    \includegraphics[width=0.45\textwidth]{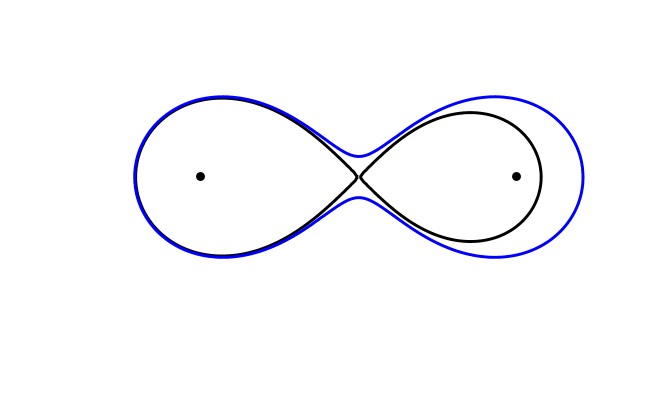} 
    \caption{Contour $\Sigma$ and $|(1-z^2)^{20+\frac{1}{2}}|=\sqrt{2}$.}
    \label{fig: curves zeros}
\end{figure}

Indeed, \eqref{EQ Zeros prev} implies
\begin{equation}\label{EQ Zeros prev2}
   (1-z)^{\frac{1}{2}}(1+z)^{\frac{1}{2}}(1-z^2)^m=i\sqrt{2}.
\end{equation}
Figure \ref{fig: curves zeros} shows, for $m=20$ the curve $|(1-z^2)^{m+\frac{1}{2}}|=\sqrt{2}$, where the zeros of $q(m,z)$ are located,  as well as the curve  $|1-z^2|=1$. It is claimed that this is the limiting curve for the zeros of the polynomials $P_m(z)$. See Corollary \ref{Corollary Asymp Zeros}.

\begin{proof}[\textbf{Proof of Theorem \ref{Thm Zeros out}}]
Let $m$ be large enough,  to establish the precise location of the zeros of $\pi_m(z)$ in the region $ \Omega_- \setminus u_{\varepsilon}$, for each zero $z_k^*$ of $q(m,z)$ a small box centered at $z_k^*$ is constructed, see Figure \ref{Box zeros}. Rouché's theorem is then used to show that there must be exactly one zero $z_k$ of $\pi_m(z)$ inside this box.
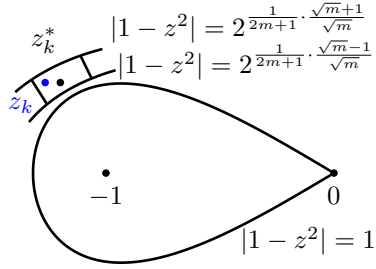
\begin{figure}[htbp]
    \centering  
    \begin{tikzpicture}[scale=2]
        \draw[line width=1pt]
        plot[domain=90:270,samples=200,smooth] ({1.4*sqrt(2)*cos(\x)},{1.4*0.3*sqrt(2)*sin(2*\x)});

        \draw[line width=1pt, domain=141:159, samples=200, smooth] 
        plot ({1.5*sqrt(1.5 + cos(2*\x))*cos(\x)}, {sqrt(1.5 + cos(2*\x))*sin(\x)});
        
        \draw[line width=1pt, domain=143:165, samples=200, smooth] 
        plot ({1.35*sqrt(1.5 + cos(2*\x))*cos(\x)}, {0.85*sqrt(1.5 + cos(2*\x))*sin(\x)});

        \coordinate (C1) at ({1.5*sqrt(1.5 + cos(2*145))*cos(145)}, {sqrt(1.5 + cos(2*145))*sin(145)});
        \coordinate (C2) at ({1.5*sqrt(1.5 + cos(2*155))*cos(155)}, {sqrt(1.5 + cos(2*155))*sin(155)});

        \coordinate (L1) at ({1.35*sqrt(1.5 + cos(2*148))*cos(148)}, {0.85*sqrt(1.5 + cos(2*148))*sin(148)});
        
        \coordinate (L2) at ({1.35*sqrt(1.5 + cos(2*159))*cos(159)}, {0.85*sqrt(1.5 + cos(2*159))*sin(159)});

        \draw[line width=1pt] (C1) -- (L1);
        \draw[line width=1pt] (C2) -- (L2);

        \node at (-1.5,0) [circle,fill,inner sep=1pt,label=below:$-1$] {};
        \node at (0,0) [circle,fill,inner sep=1pt,label=below:$0$] {};

        \node at (-1.8,0.6) [circle,fill,inner sep=1pt] {};        
        \node at (-1.7,0.88) [label=left:$z_k^*$] {};

     \node at (-1.9,0.6) [circle,fill=blue!100!,inner sep=1pt] {}; 
     
     \node at (-1.85,0.45) [label=left:{\color{blue!100!} $z_k$}] {};

    \node at (0.35,1) [label=left: \text{$|1-z^2|=2^{\frac{1}{2m+1}\cdot\frac{\sqrt{m}+1}{\sqrt{m}}}$}] {};

    \node at (0.4,0.75) [label=left: \text{$|1-z^2|=2^{\frac{1}{2m+1}\cdot\frac{\sqrt{m}-1}{\sqrt{m}}}$}] {};

         \node at (0.4,-0.45) [label=left: \text{$|1-z^2|=1$}] {};
        
    \end{tikzpicture}
    \caption{Box for the zero $z_k^*$ of $q(m,z)$.}
    \label{Box zeros}
\end{figure}

To achieve this consider the parametrization
\begin{equation}\label{EQ param 1-z^2}
 1 - z^2 = \rho e^{i\theta}, \hspace{1cm} 0 \le \theta < 2\pi.
\end{equation}
Then, from \eqref{EQ Zeros prev2} and \eqref{EQ param 1-z^2}, the zeros $z_k^*$ of $q(m,z)$ correspond to $\rho =2^{\frac{1}{2m+1}}$ and can be enumerated according to $\textnormal{Arg}(1-(z_k^*)^2)$, where
    $$\left(m+\frac{1}{2}\right)\theta_k=\frac{\pi}{2}+2\pi k.$$

For each zero $z_k^*\in  \Omega_- \setminus u_{\varepsilon}$ of $q(m,z)$, the box enclosing it is constructed as follows: the lower side of the box corresponds to a small segment of the curve
\begin{equation}\label{Eq. Lopwer box}
  |1-z^2|=2^{\frac{1}{2m+1}\cdot\frac{\sqrt{m}-1}{\sqrt{m}}} . 
\end{equation}
The upper side of the box corresponds to a small segment of the curve
\begin{equation}\label{Eq. upper box}
    |1-z^2|=2^{\frac{1}{2m+1}\cdot\frac{\sqrt{m}+1}{\sqrt{m}}}.
\end{equation}
Finally, let $\eta=\frac{\ln2}{\sqrt{m}}$, the left and right sides of the box correspond to fixing the angle $\theta$ at $$\widetilde{\theta}=\theta_k\pm\frac{\eta}{m+\frac{1}{2}},$$ respectively, and varying the parameter $\rho$ between $2^{\frac{1}{2m+1}\cdot\frac{\sqrt{m}-1}{\sqrt{m}}}$ and $2^{\frac{1}{2m+1}\cdot\frac{\sqrt{m}+1}{\sqrt{m}}}$.

Note that for $z$ on the lower side of the box, \eqref{Eq. Lopwer box} gives
\begin{equation}\label{Eq: bound q lower side box}
   |q(m,z)|=|1-\sqrt{2}i(1-z^2)^{-m-\frac{1}{2}}|\ge 1-e^{-\frac{\ln2}{2\sqrt{m}}}>\frac{\ln2}{2\sqrt{m}}.
\end{equation}
Similarly, on the upper side of the box, \eqref{Eq. upper box} gives
\begin{equation} \label{Eq: bound q upper side box}
    |q(m,z)|=|1-\sqrt{2}i\left(1-z^2\right)^{-m-\frac{1}{2}}|
    \ge e^{\frac{\ln2}{2\sqrt{m}}}-1>\frac{\ln2}{2\sqrt{m}}.
\end{equation}
On the left and right sides of the box
\begin{align*}
    |q(m,z)|&=|1-\sqrt{2}i(1-z^2)^{-m-\frac{1}{2}}|\\
    &=|1-\sqrt{2}i\rho^{-m-\frac{1}{2}}e^{-i\left(m+\frac{1}{2}\right)\theta_k\pm i\eta}|\\
    &=|1-\sqrt{2}\rho^{-m-\frac{1}{2}}e^{\pm i\eta}|\\
    &\ge |\sqrt{2}\rho^{-m-\frac{1}{2}}\sin(\pm \eta)|\\
    &\ge \frac{\sin\eta}{2^{\frac{1}{2\sqrt{m}}}}\\
    &> \frac{\sin\eta}{2}\\
    &>\frac{\ln2}{2\sqrt{m}}.
\end{align*}
Therefore $$q(m,z)> \frac{\ln2}{2\sqrt{m}}=:q_0(m),$$ on the boundary of the box. 

Now, $r(m,z)=\mathcal{O}(1/m)$, so there exist $M\in \mathbb{N}$ and a constant $C>0$, so that 
$$|r(m,z)|\le \frac{C}{m},\hspace{0.3cm}\text{for all $m\ge M$.}$$
Pick $M$ so that for $m \ge M$, $r(m,s) < q_0(m)$. Rouché's theorem implies that $q(m,z)$ and $q(m,z) + r(m,z)$ have the same number of zeros inside this box, that is, exactly one. Expression \eqref{EQ pi_n(z)=F(q+r)} shows that $\pi_m(z)$ must also have exactly one zero inside this box. One checks that the size of the box is $\mathcal{O}(1/m)$, which completes the proof of the Theorem.

\end{proof}
\subsection{Zeros near the origin}
\begin{proof}[\textbf{Proof of Theorem \ref{Thm Asymptotics local}}]
Start by obtaining a more explicit expression for the entries in the first column of the matrix $B(z)$ in \eqref{B(z) near origin} ($z\in \Omega_-\cap u_{\varepsilon}\cap \mathbb{C}^+$). From \eqref{B(z) near origin},
\begin{align*}
    B(z) &= \frac{\sqrt{2}}{2}\, N(z) \widehat{w}(z)^{\frac{\sigma_3}{2}} a^{\sigma_3}\begin{pmatrix}
            1 &1\\
            1&-1
        \end{pmatrix}\left(\xi \, e^{-\frac{5\pi i}{4}}\right)^{\frac{\sigma_3}{2}}Z\left(\xi\,e^{\frac{-5\pi i}{4}}\right)a^{-\sigma_3}e^{\frac{i\xi^2}{4}\sigma_3}\widehat{w}(z)^{-\frac{\sigma_3}{2}}V(z)\\
        &=\frac{\sqrt{2}}{2}\left[\frac{(1+z)^{\frac{1}{2}}}{z^{\frac{1}{2}}}\widehat{w}^{\frac{1}{2}}a\right]^{\sigma_3}
        \begin{pmatrix}
            1 &1\\
            1&-1
        \end{pmatrix}\left(\xi \, e^{-\frac{5\pi i}{4}}\right)^{\frac{\sigma_3}{2}}Z\left(\xi\,e^{\frac{-5\pi i}{4}}\right)a^{-\sigma_3}e^{\frac{i\xi^2}{4}\sigma_3}\widehat{w}(z)^{-\frac{\sigma_3}{2}}V(z)\\
        &=\frac{\sqrt{2}}{2}\alpha^{\sigma_3}\begin{pmatrix}
            1 &1\\
            1&-1
        \end{pmatrix}\beta^{\sigma_3}\,Z\,\gamma^{\sigma_3}
        \begin{pmatrix}
       1& 0\\
      -\frac{1}{(1-z^2)^{m}\,\widetilde{w}(z)}  & 1
    \end{pmatrix}
    \begin{pmatrix}
        (1+z)^m & 0 \\
        0 & (1+z)^{-m}
     \end{pmatrix}\\
     &=\frac{\sqrt{2}}{2}
     \begin{pmatrix}
         \alpha\beta & \alpha\beta^{-1}\\
         \alpha^{-1}\beta& -\alpha^{-1}\beta^{-1} 
     \end{pmatrix} Z 
     \begin{pmatrix}
       \gamma& 0\\
      -\frac{\gamma^{-1}}{(1-z^2)^{m}\,\widetilde{w}(z)}  & \gamma^{-1}
    \end{pmatrix}
    \begin{pmatrix}
        (1+z)^m & 0 \\
        0 & (1+z)^{-m}
     \end{pmatrix},
\end{align*}
with
\begin{equation}\label{Eq def alpha, beta, gamma}
    \alpha= \frac{(1+z)^{\frac{1}{2}}}{z^{\frac{1}{2}}}\widehat{w}^{\frac{1}{2}}a, \hspace{0.5cm} \beta=\left(\xi \, e^{-\frac{5\pi i}{4}}\right)^{\frac{1}{2}},\hspace{0.5cm}\text{and}\hspace{0,3cm}\gamma=a^{-1}e^{\frac{i\xi^2}{4}}\widehat{w}(z)^{-\frac{1}{2}}.
\end{equation}
This gives
\begin{equation*}
    B_{11}(z)=\frac{\sqrt{2}}{2}(1+z)^m\left[\gamma\left(\alpha\beta Z_{11}+\alpha\beta^{-1}Z_{21}\right)-\frac{1}{\gamma(1-z^2)^m\widetilde{w}}\left(\alpha\beta Z_{12}+\alpha\beta^{-1}Z_{22}\right)\right]
\end{equation*}
and 
\begin{equation*}
    B_{21}(z)=\frac{\sqrt{2}}{2}(1+z)^m\left[\gamma\left(\alpha^{-1}\beta Z_{11}-\alpha^{-1}\beta^{-1}Z_{21}\right)-\frac{1}{\gamma(1-z^2)^m\widetilde{w}}\left(\alpha^{-1}\beta Z_{12}-\alpha^{-1}\beta^{-1}Z_{22}\right)\right].
\end{equation*}
Using the fact that $\gamma^2(1-z^2)^m\widetilde{w}=a^{-2}$ (see \eqref{Eq def alpha, beta, gamma}), this reduces to
\begin{equation*}
    B_{11}(z)=\frac{\sqrt{2}}{2}(1+z)^m\gamma\,\alpha\,\beta\left[\left( Z_{11}+\beta^{-2}Z_{21}\right)-a^2\left( Z_{12}+\beta^{-2}Z_{22}\right)\right]
\end{equation*}
and 
\begin{equation*}
    B_{21}(z)=\frac{\sqrt{2}}{2}(1+z)^m\gamma\,\alpha\,\beta\left[\left(\alpha^{-2} Z_{11}-\alpha^{-2}\beta^{-2}Z_{21}\right)-a^2\left(\alpha^{-2}Z_{12}-\alpha^{-2}\beta^{-2}Z_{22}\right)\right].
\end{equation*}

Now, recall from \eqref{Eq Asymptotics for pi_m}:
\begin{align*}  
\pi_m(z)&=B_{11}\left[\left(1+\frac{\mathcal{E}_{11}^{(1)}}{z}\right)\widetilde{E}_{11}+ \frac{\mathcal{E}_{12}^{(1)}}{z}\widetilde{E}_{21}\right]+ B_{21}\left[\left(1+\frac{\mathcal{E}_{11}^{(1)}}{z}\right)\widetilde{E}_{12}+ \frac{\mathcal{E}_{12}^{(1)}}{z}\widetilde{E}_{22}\right],\\
&=\frac{\sqrt{2}}{2}(1+z)^m\gamma\,\alpha\,\beta\Bigg\{ \left[\Big( Z_{11} +\beta^{-2}Z_{21}\Big) - a^2\left(Z_{12}+\beta^{-2}Z_{22}\right)\right]\\
&\hspace{6cm}\cdot\left[\left(1+\frac{\mathcal{E}_{11}^{(1)}}{z}\right)\widetilde{E}_{11}+ \frac{\mathcal{E}_{12}^{(1)}}{z}\widetilde{E}_{21}\right]\\
&\hspace{1cm}+\left[\Big(\alpha^{-2}\ Z_{11} -\alpha^{-2}\beta^{-2}Z_{21}\Big) - a^2\left(\alpha^{-2} Z_{12}+\alpha^{-2}\beta^{-2}Z_{22}\right)\right]\\
&\hspace{6cm}\cdot\left[\left(1+\frac{\mathcal{E}_{11}^{(1)}}{z}\right)\widetilde{E}_{12}+ \frac{\mathcal{E}_{12}^{(1)}}{z}\widetilde{E}_{22}\right]\Bigg\}.
\end{align*}
Then, one can write
\begin{equation}\label{Eq. pi_m local}
   \pi_m(z) =\frac{\sqrt{2}}{2}(1+z)^m\gamma\,\alpha\,\beta\,\left(\,\tilde{q}_l(z)+ \tilde{r}_l(z)\right),
\end{equation}
where
\begin{equation}\label{Eq. q_l}
\begin{split}
  \tilde{q}_l(z)=& \Big(Z_{11} + \beta^{-2}Z_{21}\Big) - a^2\Big( Z_{12} + \beta^{-2}Z_{22}\Big) \\
     &\hspace{0.5cm}  - \frac{a^2}{z} \Big[\Big(\alpha^{-2} Z_{11} - \alpha^{-2}\beta^{-2}Z_{21}\Big) - a^2\Big(\alpha^{-2} Z_{12} - \alpha^{-2}\beta^{-2}Z_{22}\Big)\Big],
\end{split}
\end{equation}
and 
\begin{equation}\label{Eq. r_l}
\begin{split}
    \tilde{r}_l(z)&=\frac{\sqrt{2}}{(1+z)^m\,\gamma\,\alpha\,\beta}\Bigg\{\frac{B_{11}}{z}\left(\mathcal{E}_{11}^{(1)}\widetilde{E}_{11}+\mathcal{E}_{12}^{(1)}\widetilde{E}_{21}\right)\\
    &+B_{21}\left[\widetilde{E}_{12}+\frac{1}{z}\left(\mathcal{E}_{11}^{(1)}\widetilde{E}_{12}+(\mathcal{E}_{12}^{(1)}+a^2)\widetilde{E}_{22}\right)\right]+B_{11}(\widetilde{E}_{11}-1)-\frac{B_{21}}{z}a^2(\widetilde{E}_{22}-1)\Bigg\}.
    \end{split}
\end{equation}
Let $R>0$ be arbitrary, for $|s|=2R$, the fact that $\widetilde{E}_{11}, \widetilde{E}_{22} = 1+\mathcal{O}\left(\frac{1}{m}\right)$, $\mathcal{E}_{12}^{(1)}=-a^2+\mathcal{O}\left(\frac{1}{m}\right)$ and $\widetilde{E}_{12}, \widetilde{E}_{21},\mathcal{E}_{11}^{(1)} = \mathcal{O}\left(\frac{1}{m}\right)$ implies the existence of a constant $C>0$ and an integer $M>0$, such that if $m>M$, then
\begin{equation}\label{Eq. Bound r_l}
   |r_l(s)|\le \frac{C}{\sqrt{m}}.
\end{equation}
Hence, for $|\xi|<R$, the Cauchy integral formula gives
\begin{align*}
   \abs{r_l(\xi)}&=\left|\frac{1}{2\pi i}\oint_{|s|=2M}\frac{r_l(s)}{s-\xi}\, ds\right|\\
   &\le\frac{1}{2\pi R}\oint|r_l(s)|\,|ds|\\
   &\le\frac{1}{2\pi R}\frac{ 2\pi\, R\,C}{\sqrt{m}}\\
   &=\frac{C}{\sqrt{m}}.
\end{align*}
This show that for  $|\xi|< R$ and $z=z(\xi)\in \Omega_-\cap u_{\varepsilon}\cap \mathbb{C}^+$, $\tilde{r}_l(z)=\mathcal{O}\left(\frac{1}{\sqrt{m}}\right)$.

Now, the term $\tilde{q}_l(z)$ in \eqref{Eq. pi_m local} simplifies to
\begin{equation*}
   \tilde{q}_l(z)= \left(1-\frac{a^2}{z}\alpha^{-2}\right)(Z_{11}-a^2Z_{12})  +\beta^{-2} \left(1+\frac{a^2}{z}\alpha^{-2}\right)(Z_{21}-a^2 Z_{22}).
\end{equation*}
Using the expressions for $\alpha$, $\beta$ and $\gamma$ in \eqref{Eq def alpha, beta, gamma}, one observes that
\begin{equation*}
    1-\frac{a^2}{z}\alpha^{-2}=2+\mathcal{O}\left(z\right), \hspace{0.5cm} 1+\frac{a^2}{z}\alpha^{-2}=\mathcal{O}\left(z\right),\hspace{0.5cm}\text{and} \hspace{0.5cm}\beta^{-2}=\mathcal{O}\left(\frac{1}{\sqrt{m}}\right).
\end{equation*}
Then,  \eqref{Eq. pi_m local} yields
\begin{equation}\label{Eq pi_m local 1}
    \pi_m(z) = \frac{\sqrt{2}}{2}(1+z)^m \gamma\,\alpha\,\beta\,\Bigg\{\hspace{-0.15cm}\left(1\hspace{-0.05cm}-\hspace{-0.05cm}\frac{1}{(1-z)^{\frac{1}{2}}(1+z)^{\frac{1}{2}}}\right)\hspace{-0.05cm}\left(Z_{11}-a^2Z_{12}\right)+r_{\ell}\Bigg\},
\end{equation}
where $$r_{\ell}=r_{\ell}(m,z)=\mathcal{O}\left(\frac{1}{\sqrt{m}}\right).$$
To further simplify \eqref{Eq pi_m local 1} one uses the definition of $Z$ given in Section \ref{sec-RHP PC}.  Note that the region associated to $ \Omega_-\cap u_{\varepsilon}\cap \mathbb{C}^+$ in the $z$-plane, corresponds to  $\Omega_3$ in the variable $\eta=\xi e^{-\frac{5\pi i}{4}}$. This implies that $Z=Z_3$ in this region, see \eqref{Eq sol RHP for P}. Then, using the definition of $Z_3$ in \eqref{Eq. def Z_n} one obtains
\begin{equation}\label{Eq local zeros final}
   Z_{11}-a^2Z_{12}=\frac{1}{\sqrt{2}}\left[e^{\frac{\pi i}{4}}D_{-\frac{1}{2}}\left(e^{-\frac{3 \pi i}{4}}\xi\right)-i\sqrt{2}D_{-\frac{1}{2}}\left(e^{-\frac{5 \pi i}{4}}\xi\right)\right].
\end{equation}
Also, from \eqref{Eq def alpha, beta, gamma}
\begin{equation}\label{Eq. product alpha beta gamma}
    \alpha\,\beta\,\gamma = \frac{(1+z)^{\frac{1}{2}} \, \left(\xi\,e^{-\frac{5\pi i}{4}}\right)^{\frac{1}{2}}e^{\frac{i\xi^2}{4}}}{z^{\frac{1}{2}}}.
\end{equation}
The result in Theorem \ref{Thm Zeros in} is then obtained substituting \eqref{Eq local zeros final} and \eqref{Eq def alpha, beta, gamma} in \eqref{Eq pi_m local 1}.
\end{proof}

\begin{proof}[\textbf{Proof of Theorem \ref{Thm Zeros in}}]
Consider the disk $\Omega_R=\{\xi\in\mathbb{C}:\,\,|\xi|< R\}$, where $R>0$ is arbitrary.  Since the term $1\hspace{-0.05cm}-(1-z)^{-1/2}(1+z)^{-1/2}$ do not have zeros near the origin, equations \eqref{Eq pi_m local 1} and \eqref{Eq local zeros final} show that for $\xi\in \Omega_R$ with $z=z(\xi)\in \Omega_-\cap u_{\varepsilon}\cap \mathbb{C}^+$, the zeros of $\pi_m$ are determined by the solutions of the equation 
\begin{equation*}
   q_{\ell}(\xi)+r_l=0,
\end{equation*}
with
\begin{equation}\label{EQ Zeros prev local}
    q_{\ell}(\xi):=e^{\frac{\pi i}{4}}D_{-\frac{1}{2}}\left(e^{-\frac{3 \pi i}{4}}\xi\right)-i\sqrt{2}D_{-\frac{1}{2}}\left(e^{-\frac{5 \pi i}{4}}\xi\right).
\end{equation}
Since $r_l=\mathcal{O}(1/\sqrt{m})$, there exists a constant $C>0$ and an integer $M>0$ such that if $m\ge M$ then  
\[|r_l|<\frac{C}{\sqrt{m}}.\]
Let $\{\xi_k^*\}$ be the set of all zeros of $q_{\ell}(\xi)$ in $\in \Omega_R$, with $z(\xi_k^*)\in \Omega_-\cap u_{\varepsilon}\cap \mathbb{C}^+$ (there is a finite number of them). Note that if $|\xi-\xi_k^*|<\frac{C}{m^{1/4}}$, then
$$|q_l(\xi)|\ge \frac{\widetilde{C}}{m^{1/4}},$$
for some constant $\widetilde{C}>0$ independent of $m$.

Now, for each zero $\xi_k^*$ construct disks centered at $\xi_k^*$, with radius $r_m=\frac{C}{m^{1/4}}$. Then there exists $M>0$ such that for $m\ge M$  and $\xi$ on the boundary of the disks one has
$$r_l<\frac{C}{\sqrt{m}}<\frac{\widetilde{C}}{m^{1/4}}<q_l(\xi).$$
Here, \( M\) can be chosen large enough so that for $m\ge M$, the disks do not intersect each other. Rouché's theorem implies that $q_l(\xi)$ and $q_l(\xi) + r_l(\xi)$ have the same number of zeros inside each disk, that is, exactly one. Moreover, $q_{\ell}$ cannot have a double zero since it is the solution of a second-order linear differential equation. Equation \eqref{Eq pi_m local 1} shows that $\pi_m$ must also have exactly one zero  $z_k(\xi_k)$ inside each disk.
\end{proof}

\begin{figure}[H]
    \centering
    \includegraphics[width=0.45\textwidth]{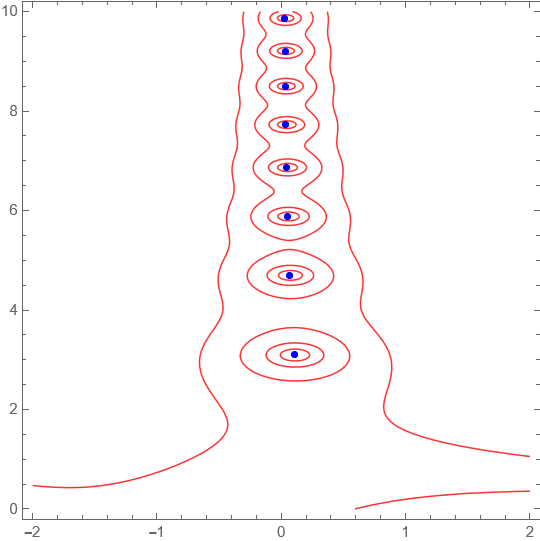} 
    \caption{Contour plot of $q_{l}(100,\xi)$ and local zeros in the $\xi$ plane.}
    \label{Temp local zeros}
\end{figure}

Unfortunately we have not been able to obtain results on the location of the zeros of the function $q_l$ to translate them into information about the zeros of the polynomials $\pi_m$, as we did in the case of the zeros away of the origin. Nonetheless, we have verified numerically that our results correctly predict the zeros of the polynomials approach the zeros of the function $q_{\ell}(m,z)$. Figure \ref{Temp local zeros} is a plot (in the $\xi$-plane, $m=100$) showing with a contour plot the location of the zeros of \eqref{EQ Zeros prev local}, along with the eight zeros of polynomial $\pi_{100}$ closest to the origin.
\section{Conclusions}
\label{Conclusions}
This work studies the asymptotic description and the zeros of the family of non-classical Jacobi polynomials arising in the works of Boros and Moll \cite{boros1999integral, boros1999sequence}. By analyzing the associated Riemann-Hilbert problem formulated in \cite{kuijlaars2005orthogonality}, explicit asymptotic formulas for the polynomials \( \pi_m(z)=\frac{1}{\kappa_{m}} P_m^{(m+1/2,-m-1/2)}(z)\) have been derived in the regions where the zeros of the polynomials are located, see Theorems \ref{Thm Asymptotics} and  \ref{Thm Asymptotics local}. These results were used to establish the limiting distribution of the zeros of \( \pi_m(z) \) as \( m \to \infty \). Besides confirming Conjecture \ref{cojecture distance to -1} and establishing that the zeros accumulate on the lemniscate \( |1 - z^2| = 1 \) in Corollary \ref{Corollary Asymp Zeros}, a complete description of the modulus and arguments of the zeros was provided in Theorem \ref{Thm Zeros out}, complementing the work of Driver and Möller \cite{driver2001zeros}. These results extend the understanding of the zeros of non-classical Jacobi polynomials.

\subsection{Future Work}
Several directions for future research arise from this work. One open question concerns the location of the zeros near the origin, which was shown in Theorem \ref{Thm Zeros in} to be determined by the solutions of the equation  
\begin{equation*}  
e^{\frac{\pi i}{4}}D_{-\frac{1}{2}}\left(e^{-\frac{3 \pi i}{4}}\xi\right)-i\sqrt{2}D_{-\frac{1}{2}}\left(e^{-\frac{5 \pi i}{4}}\xi\right)=0.  
\end{equation*}  
This result was also verified numerically in Figure \ref{Temp local zeros}, yet the methods presented in this work do not provide analytic results on the solutions of this equation. A natural question is whether a more explicit characterization of these solutions can be obtained. For example, where are the zeros of this special function equation?, and is it possible to explicitly determine the distance from the origin to the nearest zero of \( \pi_m(z) \)?

Another direction for further study involves extending the analysis of Jacobi polynomials \( P_n^{(\alpha_n, \beta_n)}(z) \) beyond the families considered by Kuijlaars et al., where the limits  
\begin{equation}\label{Eq. Limits A,B}  
\lim_{n\rightarrow \infty} \frac{\alpha_n}{n}=A\hspace{1cm}\textnormal{and} \hspace{1cm}\lim_{n\rightarrow \infty} \frac{\beta_n}{n}=B, 
\end{equation}  
are finite. Their work established that the zeros of these polynomials accumulate along curves determined by the position of \( (A,B) \) within five distinct regions of the Cartesian plane.  

\begin{figure}[H]
    \centering
    \includegraphics[width=0.4\textwidth]{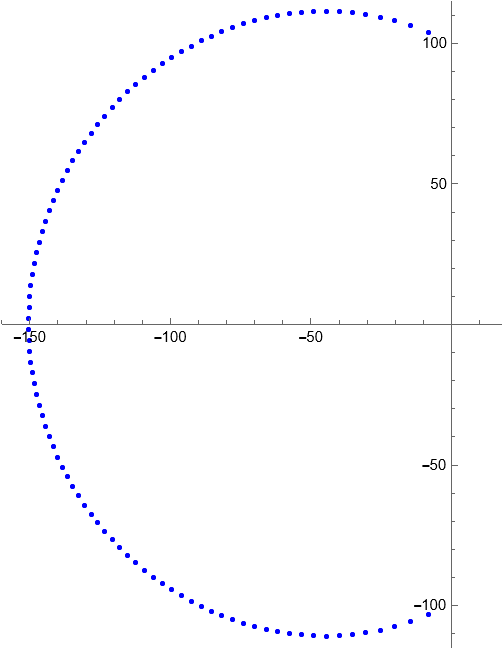} 
    \caption{zeros of $P_{100}^{\left(10000,-10000\right)}(z)$.}
    \label{fig: Zeros m square m 100}
\end{figure}

A natural question arises when these limits do not exist. For instance, if \( \alpha_n \) and \( \beta_n \) grow at order \( n^2 \), causing \( A \) and \( B \) to diverge to \( \pm \infty \). Numerical evidence suggests that, in such cases, the zeros tend to infinity while remaining confined to structured curves, which should converge to a fixed limiting shape after appropriate scaling—resembling the behavior of partial sums of the exponential function, see \cite{kriecherbauer2007locating}. A first step in exploring this phenomenon could be to analyze the case \( \alpha_m = m^2 \) and \( \beta_m = -m^2 \). The distribution of zeros for \( m = 100 \) is shown in Figure \ref{fig: Zeros m square m 100}.  




\end{document}